\documentclass[11pt, reqno]{amsart}
\usepackage{amsfonts,latexsym}
\usepackage{amsmath}
\usepackage{amscd}
\usepackage{float,amsmath,amssymb,mathrsfs,bm,multirow,graphics}
\usepackage[dvips]{graphicx}
\usepackage[percent]{overpic}
\usepackage[pdftex]{color}

\newcommand{\nequation}{\setcounter{equation}{0}}
\renewcommand{\theequation}{\mbox{\arabic{section}.\arabic{equation}}}
\newcommand{\R}{{\Bbb R}}

\newcommand{\C}{{\Bbb C}}

\newcommand{\proofbegin}{\noindent{\it Proof.\,\,}}
\newcommand{\proofend}{\hfill$\Box$\bigskip}

\newcommand{\lot}{\text{\upshape lower order terms}}

\def\XXint#1#2#3{{\setbox0=\hbox{$#1{#2#3}{\int}$}
\vcenter{\hbox{$#2#3$}}\kern-.5\wd0}}

\newtheorem{theorem}{Theorem}[section]

\newtheorem{lemma}[theorem]{Lemma}

\newtheorem{remark}[theorem]{Remark}

\newtheorem{figuretext}{Figure}

\input epsf
\title[The Unified Method III]{\sc The Unified Method: III~Non-linearizable Problems on the Interval}

\author{J. Lenells}
\address{J.L.: Department of Mathematics, Baylor University, One Bear Place \#97328, Waco, TX 76798, USA.}
\email{Jonatan\_Lenells@baylor.edu}

\author{A. S. Fokas}
\address{A.S.F.: Department of Applied Mathematics and Theoretical Physics, University of Cambridge, Cambridge CB3 0WA, United Kingdom, and Research Center of Mathematics, Academy of Athens, 11527, Greece.}
\email{T.Fokas@damtp.cam.ac.uk} 
    
\begin{document}
\begin{abstract} 
\noindent
Boundary value problems for integrable nonlinear evolution PDEs formulated on the finite interval can be analyzed by the unified method introduced by one of the authors and used extensively in the literature. The implementation of this general method to this particular class of problems yields the solution in terms of the unique solution of a matrix Riemann-Hilbert problem formulated in the complex $k$-plane (the Fourier plane), which has a jump matrix with explicit $(x,t)$-dependence involving six scalar functions of $k$, called spectral functions. Two of these functions depend on the initial data, whereas the other four depend on all boundary values. The most difficult step of the new method is the characterization of the latter four spectral functions in terms of the given initial and boundary data, i.e. the elimination of the unknown boundary values. Here, we present an effective characterization of the spectral functions in terms of the given initial and boundary data. We present two different characterizations of this problem. One is based on the analysis of the so-called global relation, on the analysis of the equations obtained from the global relation via certain transformations leaving the dispersion relation of the associated linearized PDE invariant, and on the computation of the large $k$ asymptotics of the eigenfunctions defining the relevant spectral functions.
The other is based on the analysis of the global relation and on the introduction of the so-called Gelfand-Levitan-Marchenko representations of the eigenfunctions defining the relevant spectral functions. We also show that these two different characterizations are equivalent and that in the limit when the length of the interval tends to infinity, the relevant formulas reduce to the analogous formulas obtained recently for the case of boundary value problems formulated on the half-line.
\end{abstract}
\maketitle



\tableofcontents

\section{Introduction}\nequation
This is the third in a series of papers addressing the most difficult problem in the analysis of integrable nonlinear evolution PDEs, namely the problem of expressing the so-called {\it spectral functions} in terms of the given initial and boundary conditions. In \cite{FLnonlinearizable} this problem was analyzed for the case of the half-line. In \cite{LFtperiodic} the same problem was also analyzed with the aid of the so-called Gelfand-Levitan-Marchenko (GLM) representations. Here we analyze this problem for the case of the finite interval. 

We refer the interested reader to \cite{FLnonlinearizable} for an introduction to the unified method of \cite{F1997, F2000, Fbook} and for a discussion of the difference between linearizable versus non-linearizable boundary value problems. Here, we only note that the unified method expresses the solution $q(x,t)$ of an integrable evolution PDE in terms of an integral formulated in the complex $k$-plane. This representation is similar to the integral representation obtained by the new method for the linearized version of the given nonlinear PDE, but it also contains the entries of a certain matrix-valued function $M(x,t,k)$, which is the solution of a matrix Riemann-Hilbert (RH) problem. The main advantage of the new method is the fact that this RH problem involves a jump matrix with {\it explicit} $(x,t)$-dependence, uniquely defined in terms of the spectral functions. For the problem on the interval, there exist six spectral functions denoted by $\{a(k), b(k), A(k), B(k), \mathcal{A}(k), \mathcal{B}(k)\}$. The functions $\{a(k), b(k)\}$ are defined in terms of the initial data $q_0(x) = q(x,0)$ via a system of linear Volterra integral equations; the functions $\{A(k), B(k)\}$ and $\{\mathcal{A}(k), \mathcal{B}(k)\}$ are defined in terms of the boundary values at $x = 0$ and $x = L$ respectively, also via systems of linear Volterra integral equations. However, the integral equations defining $\{A(k), B(k)\}$ and $\{\mathcal{A}(k), \mathcal{B}(k)\}$ involve {\it all} boundary values, whereas for a well-posed problem only a subset of these boundary values can be prescribed as boundary conditions.
Thus, the complete solution of a concrete initial-boundary value problem requires the characterization of $\{A(k), B(k), \mathcal{A}(k), \mathcal{B}(k)\}$ in terms of the given initial and boundary conditions. For example, for the Dirichlet problem of the NLS on the interval, $0 < x < L$, it is necessary to characterize $\{A(k), B(k), \mathcal{A}(k), \mathcal{B}(k)\}$ in terms of $q_0(x)$, $g_0(t) := q(0,t)$, and $h_0(t) = q(L,t)$.

A characterization of the spectral functions is called {\it effective} if it fulfills the following requirements: $(a)$ In the linear limit, it yields an effective solution of the linearized boundary value problem, i.e. it yields a solution in the form of an integral which involves the transforms of the given initial and boundary conditions. $(b)$ For `small' boundary conditions, it yields an effective perturbative scheme, i.e. it yields an expression in which each term can be computed uniquely in a well-defined recursive scheme.

The effective characterization presented here is based on the construction of the generalized Dirichlet to Neumann map, i.e. on the characterization of the unknown boundary values in terms of the given initial and boundary conditions. This characterization employs the same three ingredients introduced in \cite{FLnonlinearizable} for the analogous problem on the half-line: $(a)$ The computation of the large $k$ asymptotics of the eigenfunctions $\Phi(t,k)$ and $\varphi(t,k)$ defining $\{A(k), B(k)\}$ and $\{\mathcal{A}(k), \mathcal{B}(k)\}$ respectively. $(b)$ The analysis of the so-called {\it global relation} and of the equations obtained from the global relation under the transformations which leave invariant the dispersion relation of the associated linearized equation.
$(c)$ The construction of a perturbative scheme for establishing effectiveness. 

This paper is organized as follows: In section \ref{prelsec} we review the main results of \cite{FI2004}. In section \ref{DirichletNeumannsec}, by employing the ingredients $(a)$-$(c)$ mentioned earlier, we analyze both the Dirichlet and the Neumann problems of the NLS on the finite interval with zero initial conditions. In the former case we express $\{g_1(t) = q_x(0,t), h_1(t) =q_x(L,t)\}$ in terms of $\{g_0(t), h_0(t)\}$, whereas in the latter case we express $\{g_0(t), h_0(t)\}$ in terms of $\{g_1(t), h_1(t)\}$. Furthermore, in section \ref{DirichletNeumannsec} we show that as $L \to \infty$, the formulas for $g_1$ and $g_0$ respectively, coincide with the analogous formulas obtained in \cite{FLnonlinearizable}.
In section \ref{glmsec} we analyze the same problems, but we now express the unknown boundary values in terms of the GLM representations. In this connection we correct the expressions for $\{g_1(t), h_1(t)\}$ in terms of the GLM representations presented in \cite{FI2004}. Furthermore, in section \ref{glmsec} we establish the equivalence of the formulas obtained directly (i.e. the formulas obtained in section \ref{DirichletNeumannsec}) with the formulas obtained via the GLM representations (i.e. the formulas presented in section \ref{glmsec}). 

We emphasize that for problems on the finite interval it is necessary to employ {\it all} three ingredients mentioned earlier. In particular, the asymptotics of $\Phi(t,k)$ and $\varphi(t,k)$ yield {\it several} possible formulas for the unknown boundary values, so it is absolutely necessary to employ the ingredient $(c)$ in order to choose the one that yields an effective solution. As a warning, we present in appendix \ref{noneffectiveapp} a particular set of formulas for $\{g_1(t), h_1(t)\}$ which do {\it not} yield an effective solution.

\section{Preliminaries}\nequation\label{prelsec}
We consider the NLS equation on the interval $[0,L]$:
\begin{align}\label{nls}
 & iq_t + q_{xx} - 2 \lambda |q|^2 q = 0, \qquad \lambda = \pm 1, \; 0 < x < L, \; 0 < t < T,
\end{align}
where $L>0$ is the length of the interval and $T>0$ is a fixed finite time.
We let $g_0(t)$ and $h_0(t)$ denote the Dirichlet boundary values of $q(x,t)$, whereas $g_1(t)$ and $h_1(t)$ denote the Neumann boundary values:
\begin{align}\label{g0h0g1h1def}
q(0,t) = g_0(t), \quad q(L,t) = h_0(t), \quad q_x(0,t) = g_1(t), \quad q_x(L,t) = h_1(t).
\end{align}

\subsection{Bounded and analytic eigenfunctions}
In what follows we present a summary of the results obtained in \cite{FI2004}.

\begin{figure}
\begin{center}
\medskip
 \begin{overpic}[width=.3\textwidth]{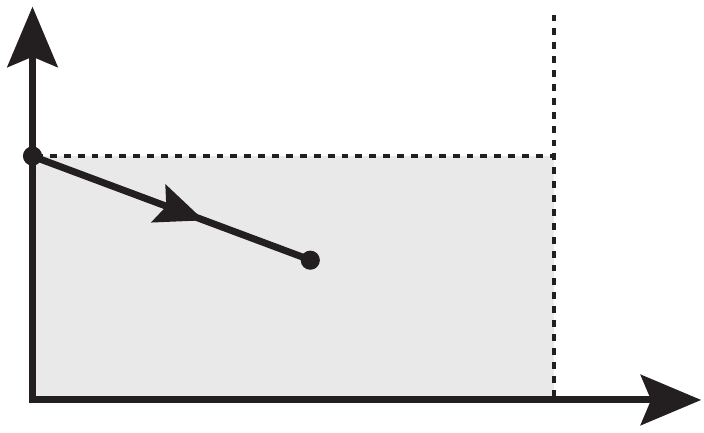}
 \put(-27,36){\small $(x_1,t_1)$}
 \put(47,17){\small $(x,t)$}
 \put(38,-8){\small (1)}
   \end{overpic}
   \hspace{2.5cm}
  \begin{overpic}[width=.3\textwidth]{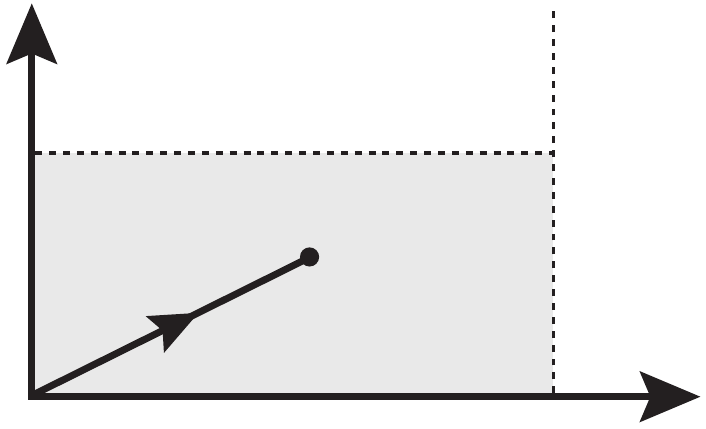}
    \put(-25,-1){\small $(x_2, t_2)$}
    \put(46,18){\small $(x,t)$} 
 \put(38,-8){\small (2)}
  \end{overpic} 
\\
\vspace{.8cm}\hspace{.08cm}
  \begin{overpic}[width=.3\textwidth]{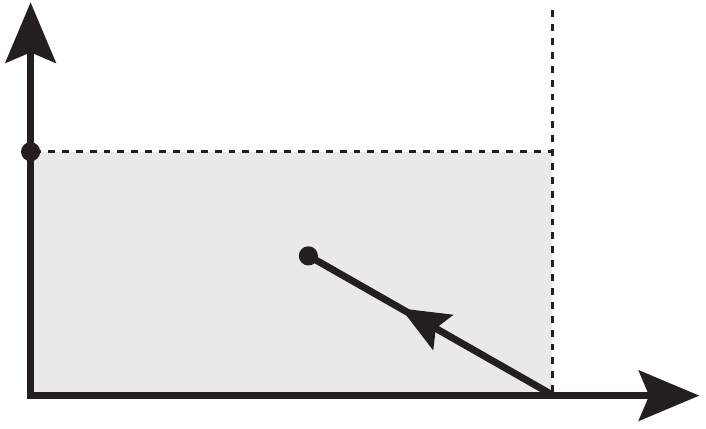}
     \put(62,-6){\small $(x_3, t_3)$}
    \put(21,20){\small $(x,t)$} 
 \put(38,-8){\small (3)}
  \end{overpic}
  \hspace{2.5cm}
    \begin{overpic}[width=.3\textwidth]{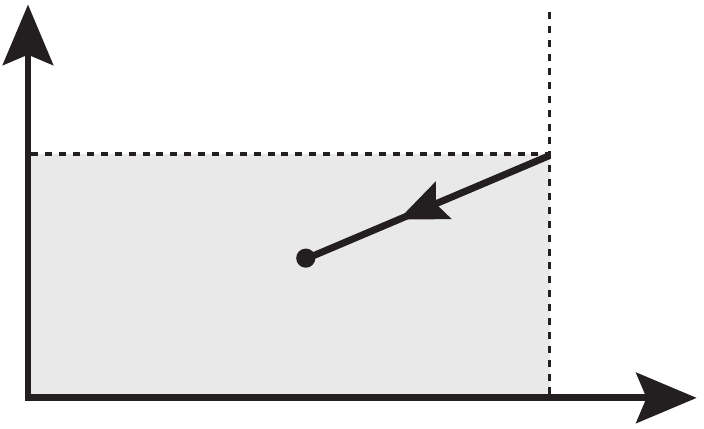}
     \put(81,36){\small $(x_4, t_4)$}
    \put(21,18){\small $(x,t)$} 
 \put(38,-8){\small (4)}
  \end{overpic}
  \bigskip
   \begin{figuretext}\label{mucontoursfig}
     The contours used to define $\{\mu_j\}_1^4$.
   \end{figuretext}
   \end{center}
\end{figure}

The Lax pair of (\ref{nls}) can be written in differential form as 
\begin{align}\label{laxdiffform}  
  d\left(e^{i(kx + 2k^2 t)\hat{\sigma}_3} \mu(x,t,k)\right) = W(x,t,k),
\end{align}
where $k \in \C$ is the spectral parameter, the closed $1$-form $W$ is defined by
\begin{align}\label{Wdef}
 & W(x,t,k) = e^{i(kx + 2k^2 t)\hat{\sigma}_3}(Q dx + \tilde{Q} dt)\mu, \qquad \sigma_3 = \begin{pmatrix} 0 & 1 \\ 1 & 0 \end{pmatrix},
  	\\ \nonumber
&    Q = \begin{pmatrix} 0 & q \\ \lambda \bar{q} & 0 \end{pmatrix}, \qquad
   \tilde{Q} = 2kQ -i(Q_x + \lambda |q|^2 )\sigma_3,
\end{align}
and $\hat{\sigma}_3A = [\sigma_3, A]$. We define four eigenfunctions $\{\mu_j\}_1^4$ by integrating from the four corners of the domain $\{0 < x < L, \; 0 < t < T\}$:
\begin{align}\label{mujdef}
\mu_j(x,t,k) = I + \int_{(x_j, t_j)}^{(x,t)} e^{-i(kx + 2k^2t)\hat{\sigma}_3} W_j(x',t',k),
\end{align}
where $W_j$ is the differential form defined in (\ref{Wdef}) with $\mu$ replaced by $\mu_j$ and $\{(x_j, t_j)\}_1^4$ denote the points $(0,T)$, $(0,0)$, $(L,0)$, and $(L,T)$, respectively, see figure \ref{mucontoursfig}. 
The spectral functions $\{a(k), b(k), A(k), B(k), \mathcal{A}(k), \mathcal{B}(k)\}$ are defined for $k \in \C$ by
\begin{align*}
& s(k) = \begin{pmatrix}\overline{a(\bar{k})} & b(k) \\
  \lambda \overline{b(\bar{k})} & a(k) \end{pmatrix}, \quad  S(k) = \begin{pmatrix}\overline{A(\bar{k})} & B(k) \\
  \lambda \overline{B(\bar{k})} & A(k) \end{pmatrix},
\quad  S_L(k) = \begin{pmatrix}\overline{\mathcal{A}(\bar{k})} & \mathcal{B}(k) \\
  \lambda \overline{\mathcal{B}(\bar{k})} & \mathcal{A}(k) \end{pmatrix}, 
\end{align*}  
where
\begin{align*}
s(k) := \mu_3(0,0,k), \qquad S(k) := \mu_1(0,0,k), \qquad S_L(k) := \mu_4(L,0,k).
\end{align*}
We also introduce the functions $\{\Phi_1, \Phi_2, \varphi_1, \varphi_2\}$ as follows:
\begin{align}\label{Phivarphidef}
& \mu_2(0,t,k) = \begin{pmatrix} \overline{\Phi_2(t, \bar{k})} & \Phi_1(t,k) \\
\lambda \overline{\Phi_1(t, \bar{k})} & \Phi_2(t,k) 
\end{pmatrix}, \qquad
\mu_3(L,t,k) = \begin{pmatrix} \overline{\varphi_2(t, \bar{k})} & \varphi_1(t,k) \\
\lambda \overline{\varphi_1(t, \bar{k})} & \varphi_2(t,k) 
\end{pmatrix}.
\end{align}

Let $\{D_j\}_1^4$ denote the four quadrants of the complex $k$-plane. 
Then,
\begin{itemize}
\item $a(k)$ and $b(k)$ are entire functions of $k$ which are bounded in $D_1 \cup D_2$.
\item $A(k)$, $B(k)$, $\mathcal{A}(k)$, $\mathcal{B}(k)$ are entire functions of $k$ which are bounded for $k \in D_1 \cup D_3$. 
\item $\Phi_1(t, k)$, $\Phi_2(t, k)$, $\varphi_1(t,k)$, $\varphi_2(t,k)$ are entire functions of $k$ which are bounded forÊ $k \in D_2 \cup D_4$. 
\end{itemize}

\subsection{The global relation}
The eigenfunctions $\{\Phi_j, \varphi_j\}_1^2$ satisfy the following {\it global relation}:
\begin{align}\label{GR}
c(t,k) = &\; \overline{\Phi_2(t,\bar{k})}\left[b(k)\overline{\varphi_2(t, \bar{k})} e^{-4ik^2t} -\overline{a(\bar{k})} \varphi_1(t,k) e^{2ikL}\right]
	\\\nonumber
& + \Phi_1(t,k) \left[a(k)\overline{\varphi_2(t, \bar{k})} - \lambda \overline{b(\bar{k})} \varphi_1(t,k) e^{4ik^2t + 2ikL}\right], \qquad k \in \C,
\end{align}
where $c(t,k)$ is an entire function of $k$ such that
\begin{align}\label{cO}
  c(t,k) = O\biggl(\frac{1 + e^{2ikL}}{k}\biggr), \qquad k \to \infty, \; k \in \C.
\end{align}  
Indeed, for each $t \in (0,T)$, let $R(x,t,k)$ be the solution of the $x$-part of the Lax pair of (\ref{nls}) such that $R(L,t,k) = I$, i.e. $R$ is the unique solution of the Volterra integral equation
\begin{align}\label{Rdef}
R(x,t,k) = I + \int_L^x e^{ik(x' - x)\hat{\sigma}_3}(QR)(x',t,k) dx, \qquad 0 < x < L.
\end{align}
It follows that $R$ is related to $\mu_3$ by
\begin{align}\label{Rmu3relation}  
  R(x,t,k) = \mu_3(x,t,k)e^{ik(L-x)\hat{\sigma}_3}\mu_3^{-1}(L,t,k).
\end{align}
Since
$$\mu_3(x,t,k) = \mu_2(x,t,k)e^{-i(kx + 2k^2t)\hat{\sigma}_3}s(k),$$
this yields
\begin{align}\label{Rmu2mu3}  
  R(x,t,k) = \mu_2(x,t,k)\bigl[e^{-i(kx + 2k^2t)\hat{\sigma}_3}s(k)\bigr]e^{ik(L-x)\hat{\sigma}_3}\mu_3^{-1}(L,t,k).
\end{align}
Evaluating the $(12)$ entry of (\ref{Rmu2mu3}) at $x = 0$ and recalling the definitions (\ref{Phivarphidef}) of $\{\Phi_j, \varphi_j\}_1^2$, we find (\ref{GR}) with $c(t,k) := R_{12}(t,k)$. It follows from (\ref{Rdef}) that $c(t,k)$ is an entire function satisfying (\ref{cO}).

In the case of vanishing initial data, the global relation (\ref{GR}) reduces to
\begin{align}\label{globalrelation}
c = \Phi_1\bar{\varphi}_2 - \bar{\Phi}_2\varphi_1 e^{2ikL}.
\end{align}

%It follows from (\ref{Rmu3relation}) that $R$ satisfies
%$$R_t = -2ik^2[\sigma_3, R] + \tilde{Q}(0,t,k)R - Re^{ikL\hat{\sigma}_3} \tilde{Q}(L,t,k).$$

\subsection{Asymptotics}
Integration by parts in (\ref{mujdef}) shows that 
%the $\mu_j$'s have the following asymptotics as $k \to \infty$:
%\begin{align}\nonumber
% & (\mu_j(x,t,k))_{12} = \frac{q}{2ik} + \biggl(\frac{q_x}{4} - \frac{i}{2}q \int_{(x_j, t_j)}^{(x,t)}\omega\biggr)\frac{1}{k^2} + O \Bigl(\frac{1}{k^3} \Bigr), \qquad k \to \infty,
%  	\\ \label{intervalmujasymptotics}
% &  (\mu_j(x,t,k))_{22} = 1 + \frac{1}{k} \int_{(x_j, t_j)}^{(x,t)}\omega + O \Bigl(\frac{1}{k^2} \Bigr), \qquad k \to \infty,
%\end{align}
%where the closed one-form $\omega$ is defined by
%$$\omega = \frac{\lambda}{2}\left[-i|q|^2dx + (\bar{q}q_x - q\bar{q}_x)dt\right],$$
%and the expansions in (\ref{intervalmujasymptotics}) are valid for $k$ approaching $\infty$ within the regions of boundedness of $(\mu_j)_{12}$ and $(\mu_j)_{22}$. Hence,
\begin{align}\nonumber
%& a(k) = 1 + \frac{a^{(1)}}{k} + O\Bigl(\frac{1}{k^2} \Bigr), \qquad k \to \infty, \quad k \in D_1 \cup D_2,
%	\\\nonumber
%& b(k) = - \frac{ig_0(0)}{2k} + \frac{1}{k^2}\left(\frac{g_1(0)}{4} - \frac{ig_0(0)}{2}\int_{(L, 0)}^{(0,0)} \omega\right) + O \Bigl(\frac{1}{k^3} \Bigr), \qquad k \to \infty, \quad k \in D_1 \cup D_2,
%	\\\nonumber
& \Phi_1(t, k) = \frac{\Phi_1^{(1)}(t)}{k} + \frac{\Phi_1^{(2)}(t)}{k^2} + O \Bigl(\frac{1}{k^3} \Bigr) + O\Bigl(\frac{e^{-4ik^2t}}{k}\Bigr), \qquad k \to \infty, \quad k \in D_2 \cup D_4,
	\\ \label{asymptoticformulas}
& \Phi_2(t, k) = 1 + \frac{\Phi_2^{(1)}(t)}{k}  + O \Bigl(\frac{1}{k^2} \Bigr), \qquad k \to \infty, \quad k \in D_2 \cup D_4,
	\\ \nonumber
& \varphi_1(t, k) = \frac{\varphi_1^{(1)}(t)}{k} + \frac{\varphi_1^{(2)}(t)}{k^2} + O \Bigl(\frac{1}{k^3} \Bigr) + O\Bigl(\frac{e^{-4ik^2t}}{k}\Bigr), \qquad k \to \infty, \quad k \in D_2 \cup D_4,
	\\ \nonumber
& \varphi_2(t, k) = 1 + \frac{\varphi_2^{(1)}(t)}{k}  + O \Bigl(\frac{1}{k^2} \Bigr), \qquad k \to \infty, \quad k \in D_2 \cup D_4,
\end{align}
where 
%$a^{(1)} =  \int_{(L, 0)}^{(0,0)} \omega$ and
\begin{align}\label{Phivarphicoefficients}
& \Phi_1^{(1)}(t) = \frac{g_0(t)}{2i}, \qquad
\Phi_1^{(2)}(t) = \frac{g_1(t)}{4} + \frac{g_0(t)}{2i}\int_{(0, 0)}^{(0,t)} \omega, \qquad
\Phi_2^{(1)}(t) = \int_{(0, 0)}^{(0,t)} \omega,
	\\  \nonumber
& \varphi_1^{(1)}(t) = \frac{h_0(t)}{2i}, \qquad
\varphi_1^{(2)}(t) = \frac{h_1(t)}{4} + \frac{h_0(t)}{2i}\int_{(L, 0)}^{(L,t)} \omega, \qquad
\varphi_2^{(1)}(t) = \int_{(L, 0)}^{(L,t)} \omega,
\end{align}
and the closed one-form $\omega$ is defined by
$$\omega = \frac{\lambda}{2}\left[-i|q|^2dx + (\bar{q}q_x - q\bar{q}_x)dt\right].$$
In particular, we find the following expressions for the boundary values:
\begin{subequations}
\begin{align}\label{g0h0Phi}
& g_0(t) = 2i \Phi_1^{(1)}(t), && h_0(t) = 2i \varphi_1^{(1)}(t),
	\\ \label{g1h1Phi}
& g_1(t) = 4 \Phi_1^{(2)}(t) + 2ig_0\Phi_2^{(1)}(t), &&
 h_1(t) = 4 \varphi_1^{(2)}(t) + 2ih_0\varphi_2^{(1)}(t), \qquad 0 < t < T.
\end{align}
\end{subequations}

We will also need the asymptotics of $c$.

\begin{lemma}\label{clemma}
The global relation (\ref{GR}) implies that 
\begin{align}\label{casymptotics} 
 c(t, k) = &\; \frac{\Phi_1^{(1)}(t)}{k} + \frac{\Phi_1^{(2)}(t) + \Phi_1^{(1)}(t)\bigl(a^{(1)} + \bar{\varphi}_2^{(1)}(t)\bigr)}{k^2} + O \Bigl(\frac{1}{k^3} \Bigr)
 	\\ \nonumber
&- \biggl[\frac{\varphi_1^{(1)}(t)}{k} + \frac{\varphi_1^{(2)}(t) + \varphi_1^{(1)}(t)\bigl(\bar{a}^{(1)} + \bar{\Phi}_2^{(1)}(t)\bigr)}{k^2} + O \Bigl(\frac{1}{k^3} \Bigr)\biggr]e^{2ikL},
%	\\ \nonumber
%& + O \Bigl(\frac{(1 + e^{2ikL})e^{4ik^2t}}{k^2} \Bigr) 
%+  O \Bigl(\frac{(1 + e^{2ikL})e^{4ik^2(T-t)}}{k} \Bigr), 
	\\ \nonumber
& \hspace{8cm} k \to \infty, \; k \in D_1 \cup D_3,
 \end{align}
 where $a^{(1)} =  \int_{(L, 0)}^{(0,0)} \omega$.
\end{lemma}
\proofbegin
See appendix \ref{clemmaapp}.
\proofend

\begin{remark}\upshape
We emphasize that the proof of lemma \ref{clemma} does {\it not} require the knowledge of the explicit form of $c(t,k)$; it only requires the \textit{existence} of an entire function $c(t,k)$ satisfying (\ref{GR}) and (\ref{cO}).
\end{remark}

\section{The Dirichlet and Neumann problems}\nequation\label{DirichletNeumannsec}
We will use the following notations:
\begin{itemize}
\item 
For $j = 1, \dots, 4$, $\partial D_j$ denotes the boundary of the $j$'th quadrant $D_j$, oriented so that $D_j$ lies to the left of $\partial D_j$. 

\item $\partial D_3^0$ denotes the contour obtained by deforming the contour $\partial D_3$ so that it passes below the zeros of $\Delta(k)$ in $\R_{\leq 0}$, i.e. below the set $\bigl\{-\frac{n\pi}{2L}\, \big| \,n = 0,1,2, \dots\bigr\}$.
Moreover, we let $\partial D_1^0 = -\partial D_3^0$.

\item The functions $f_+(k)$ and $f_-(k)$ denote the following even and odd combinations of the function $f(k)$:
\begin{align*}
 & f_+(k) = f(k) + f(-k), \qquad f_-(k) = f(k) - f(-k), \qquad k \in \C.
\end{align*}  

\item $\Delta(k)$ and $\Sigma(k)$ are defined by 
$$\Delta(k) = e^{2ikL} - e^{-2ikL}, \qquad \Sigma(k) = e^{2ikL} + e^{-2ikL}.$$ 
\end{itemize}

The following theorem expresses the spectral functions $\{A(k), B(k), \mathcal{A}(k), \mathcal{B}(k)\}$ in terms of the given boundary data via the solution of a system of nonlinear integral equations. For simplicity, we assume that $q_0(x)$ vanishes identically.

\begin{theorem}\label{th1}
Let $T < \infty$ and let $q_0(x) = 0$ for $x \geq 0$. For the Dirichlet problem it is assumed that the functions $g_0(t)$ and $h_0(t)$, $0 \leq t < T$, have sufficient smoothness and are compatible with the initial data. For the Neumann problem it is assumed that the functions $g_1(t)$ and $h_1(t)$, $0 \leq t < T$, have sufficient smoothness and are compatible with the initial data.

Then the spectral functions $\{A, B, \mathcal{A}$, $\mathcal{B}\}$ are given by
\begin{subequations}\label{ABexpressions}
\begin{align}
& A(k)  = \overline{\Phi_2(T, \bar{k})} 
\qquad B(k) = -\Phi_1(T, k)e^{4ik^2T},
	\\
& \mathcal{A}(k)  = \overline{\varphi_2(T, \bar{k})} 
\qquad \mathcal{B}(k) = -\varphi_1(T, k)e^{4ik^2T}, \qquad k \in \C,
\end{align}
\end{subequations}
where the complex-valued functions $\{\Phi_j(t, k), \varphi_j(t,k)\}_1^2$ satisfy the following system of nonlinear integral equations:
\begin{subequations}\label{Phivarphieqs}
\begin{align} 
&  \Phi_1(t, k) = \int_0^t e^{4ik^2(t' - t)} \bigl[-i\lambda |g_0|^2 \Phi_1 + (2kg_0 + ig_1)\Phi_2\bigr] (t', k) dt',
	\\
&   \Phi_2(t, k) = 1 + \lambda \int_0^t \bigl[(2k\bar{g}_0  - i\bar{g}_1)\Phi_1 + i|g_0|^2\Phi_2\bigr](t', k) dt', 
	\\ 
&  \varphi_1(t, k) = \int_0^t e^{4ik^2(t' - t)} \bigl[-i\lambda |h_0|^2 \varphi_1 + (2kh_0 + ih_1)\varphi_2\bigr] (t', k) dt',
	\\ 
&   \varphi_2(t, k) = 1 + \lambda \int_0^t \bigl[(2k\bar{h}_0  - i\bar{h}_1)\varphi_1 + i|h_0|^2\varphi_2\bigr](t', k) dt', \qquad 0 < t < T, \; k \in \C.
\end{align}
\end{subequations}

\begin{itemize}
\item[$(a)$] For the Dirichlet problem, the unknown Neumann boundary values $g_1(t)$ and $h_1(t)$ are given by the following expressions:
\begin{subequations}\label{g1h1expressions}
\begin{align}\nonumber
g_1(t) =  \frac{4}{\pi i}&\int_{\partial D_3^0} \biggl\{\frac{\Sigma(k)}{2\Delta(k)}\left[k\Phi_{1-}(t,k)  + ig_0(t) \right] - \frac{1}{\Delta(k)}\left[k\varphi_{1-}(t,k) + ih_0(t)\right]
	\\ \nonumber
& + \frac{g_0(t)}{4i}\varphi_{2-}(t,k) + \frac{ig_0(t)}{2} \Phi_{2-}(t,k) 
	\\ \label{g1expression}
& + \frac{k}{\Delta(k)} \left[(\overline{\varphi_2(t, \bar{k})} -1)\Phi_1(t, k)e^{-2ikL} - \varphi_1(t, k) (\overline{\Phi_2(t, \bar{k})} -1)\right]_- \biggr\}dk
\end{align}
and
\begin{align}\nonumber
h_1(t) = \frac{4}{\pi i}&\int_{\partial D_3^0} \biggl\{ -\frac{\Sigma(k)}{2\Delta(k)} \left[k\varphi_{1-}(t,k) + ih_0(t)\right]
+ \frac{1}{\Delta(k)}\left[k\Phi_{1-}(t,k) + ig_{0}(t)\right]
	\\ \nonumber
&  + \frac{h_0(t)}{4i}\Phi_{2-}(t,k) + \frac{ih_0(t)}{2} \varphi_{2-} (t,k)
	\\ \label{h1expression}
& + \frac{k}{\Delta(k)} \left[(\overline{\varphi_2(t, \bar{k})} -1)\Phi_1(t, k) - \varphi_1(t, k) (\overline{\Phi_2(t, \bar{k})} -1) e^{2ikL}\right]_-\biggr\}dk.
\end{align}
\end{subequations}

\item[$(b)$] For the Neumann problem, the unknown boundary values $g_0(t)$ and $h_0(t)$ are given by the following expressions:
\begin{subequations}\label{g0h0expressions}
\begin{align}\label{g0expression}
 g_0(t) =& \; \frac{2}{\pi} \int_{\partial D_3^0} \frac{1}{\Delta(k)} \biggl\{ \frac{\Sigma(k)}{2 }\Phi_{1+}(t,k) - \varphi_{1+}(t,k)
	\\ \nonumber
& + \left[\Phi_1(t,k) (\overline{\varphi_2(t,\bar{k})} -1)e^{-2ikL} - (\overline{\Phi_2(t,\bar{k})}  -1)\varphi_1(t,k) \right]_+\biggr\}dk
\end{align}
and
\begin{align} \label{h0expression}
h_0(t) = & \; \frac{2}{\pi} \int_{\partial D_3^0} \frac{1}{\Delta(k)}\biggl\{ - \frac{\Sigma(k)}{2}\varphi_{1+}(t,k)   + \Phi_{1+}(t,k)
	\\ \nonumber
&  + \left[(\overline{\varphi_2(t,\bar{k})} -1)\Phi_1(t,k) - e^{2ikL}(\overline{\Phi_2(t,\bar{k})}  -1)\varphi_1(t,k)\right]_+ \biggr\} dk.
\end{align}
\end{subequations}
\end{itemize}
\end{theorem}
\proofbegin
Equations (\ref{ABexpressions}) and (\ref{Phivarphieqs}) follow from the definitions and (\ref{mujdef}).

$(a)$ In order to derive (\ref{g1expression}) we note that the first of equations (\ref{g1h1Phi}) expresses $g_1$ in terms of $\Phi_2^{(1)}$ and $\Phi_1^{(2)}$. Furthermore, equations (\ref{asymptoticformulas}) and Cauchy's theorem imply
\begin{align}\label{Phi21fromPhi2}
  -\frac{i\pi}{2} \Phi_2^{(1)}(t) = \int_{\partial D_2} [\Phi_2 -1] dk = \int_{\partial D_4} [\Phi_2 - 1] dk
\end{align}
and
\begin{align}\label{Phi12fromPhi1}
  -\frac{i\pi}{2} \Phi_1^{(2)}(t) = \int_{\partial D_2} \bigl[k\Phi_1 - \Phi_1^{(1)} \bigr] dk 
  = \int_{\partial D_4} \bigl[k\Phi_1 - \Phi_1^{(1)}\bigr] dk.
\end{align}
Thus,
\begin{align}\label{ipiPhi21}
  i\pi \Phi_2^{(1)}(t) & = -\biggl(\int_{\partial D_2} + \int_{\partial D_4}\biggr) [\Phi_2 -1] dk
  = \biggl(\int_{\partial D_3} + \int_{\partial D_1}\biggr) [\Phi_2 -1] dk
  	\\ \nonumber
&  = \int_{\partial D_3} [\Phi_2(t,k) -1] dk - \int_{\partial D_3} [\Phi_2(t,-k) -1] dk 
    = \int_{\partial D_3} \Phi_{2-}(t,k) dk
\end{align}
and
\begin{align}\nonumber
i\pi \Phi_1^{(2)}(t) & = \biggl(\int_{\partial D_3} + \int_{\partial D_1}\biggr) \bigl[k\Phi_1 - \Phi_1^{(1)}\bigr] dk
  = \int_{\partial D_3} \bigl[k\Phi_1 - \Phi_1^{(1)}\bigr]_- dk
  	\\ \nonumber
&  = \int_{\partial D_3^0}  \biggl\{k\Phi_1 - \Phi_1^{(1)} + \frac{2e^{-2ikL}}{\Delta}\bigl[k\Phi_1 - \Phi_1^{(1)} \bigr] \biggr\}_- dk + I(t)
	\\ \label{ipiPhi12}
&  = \int_{\partial D_3^0} \biggl\{\frac{k \Sigma \Phi_1}{\Delta} - \frac{\Sigma}{\Delta} \Phi_1^{(1)} \biggr\}_- dk + I(t),
\end{align}
where the function $I(t)$ is defined by
$$I(t) = - \int_{\partial D_3^0} \biggl\{ \frac{2e^{-2ikL}}{\Delta}\bigl[k\Phi_1 - \Phi_1^{(1)}\bigr] \biggr\}_-dk.$$
The last step involves using the global relation (\ref{globalrelation}) to compute $I(t)$:
\begin{align}\nonumber
I(t) = & \int_{\partial D_3^0} \biggl\{ -\frac{2e^{-2ikL}}{\Delta}\biggl[kc - \Phi_1^{(1)} - \frac{\Phi_1^{(1)}\bar{\varphi}_2^{(1)}}{k} + \varphi_1^{(1)}e^{2ikL}\biggr] \biggr\}_- dk
	\\ \label{laststep}
& + \int_{\partial D_3^0} \biggl\{ \frac{2e^{-2ikL}}{\Delta}\frac{\Phi_1^{(1)}\bar{\varphi}_2^{(1)}}{k} -\frac{2}{\Delta}[k\varphi_1 - \varphi_1^{(1)}]\biggr\}_-dk
	\\ \nonumber
& + \int_{\partial D_3^0} \frac{2k}{\Delta}\bigl[\Phi_1 (\bar{\varphi}_2 - 1) e^{-2ikL} - (\bar{\Phi}_2 -1)\varphi_1\bigr]_-dk.
\end{align}
The asymptotics (\ref{casymptotics}) of $c(t,k)$ and Cauchy's theorem imply that the first integral on the rhs of (\ref{laststep}) equals $- i\pi \Phi_1^{(2)}(t)$. Moreover, analogously to (\ref{ipiPhi21}), we have the identity
\begin{align}\label{ipivarphi21}
    i\pi \varphi_2^{(1)}(t) = \int_{\partial D_3} \varphi_{2-}(t,k) dk.
\end{align}
Using this identity, the fact that $\bar{\varphi}_2^{(1)} = -\varphi_2^{(1)}$, and the expressions (\ref{Phivarphicoefficients}) for $\Phi_1^{(1)}$ and $\varphi_1^{(1)}$, we can write the second integral on the rhs of (\ref{laststep}) as
\begin{align*}
    \int_{\partial D_3^0} \biggl\{\frac{g_0}{2i}\varphi_{2-} -\frac{2}{\Delta}[k\varphi_{1-} + ih_0]\biggr\} dk.
\end{align*}
Therefore, equations (\ref{ipiPhi12}) and (\ref{laststep}) imply
\begin{align}\nonumber
&  2 i\pi \Phi_1^{(2)}(t) =  \int_{\partial D_3^0} \biggl\{\frac{k \Sigma \Phi_1}{\Delta} - \frac{\Sigma}{\Delta} \Phi_1^{(1)} \biggr\}_- dk
+   \int_{\partial D_3^0} \biggl\{\frac{g_0}{2i}\varphi_{2-} -\frac{2}{\Delta}[k\varphi_{1-} + ih_0]\biggr\}dk
	\\ \label{Phi12expression}
& + \int_{\partial D_3^0} \frac{2k}{\Delta}\bigl[\Phi_1 (\bar{\varphi}_2 - 1) e^{-2ikL} - (\bar{\Phi}_2 -1)\varphi_1\bigr]_-dk.
\end{align}
Using (\ref{ipiPhi21}) and (\ref{Phi12expression}) in the first of equations (\ref{g1h1Phi}), we find (\ref{g1expression}).

The expression (\ref{h1expression}) for $h_1(t)$ can be derived in a similar way. Indeed, the second of equations (\ref{g1h1Phi}) expresses $h_1$ in terms of $\varphi_2^{(1)}$ and $\varphi_1^{(2)}$. The coefficient $\varphi_2^{(1)}$ is given by (\ref{ipivarphi21}), whereas $\varphi_1^{(2)}$ satisfies the following analog of equation (\ref{ipiPhi12}):
\begin{align}\nonumber
i\pi \varphi_1^{(2)}(t) 
&  = \int_{\partial D_3} \bigl[k\varphi_1 - \varphi_1^{(1)}\bigr]_- dk
  	\\ \nonumber
&  = \int_{\partial D_3^0}  \biggl\{k\varphi_1 - \varphi_1^{(1)} - \frac{2e^{2ikL}}{\Delta}\bigl[k\varphi_1 - \varphi_1^{(1)} \bigr] \biggr\}_- dk + J(t)
	\\ \label{ipivarphi12}
&  = \int_{\partial D_3^0} \biggl\{-\frac{k \Sigma \varphi_1}{\Delta} + \frac{\Sigma}{\Delta} \varphi_1^{(1)} \biggr\}_- dk + J(t),
\end{align}
where the function $J(t)$ is defined by
$$J(t) = \int_{\partial D_3^0} \biggl\{ \frac{2e^{2ikL}}{\Delta}\bigl[k\varphi_1 - \varphi_1^{(1)}\bigr] \biggr\}_-dk.$$
The last step involves using the global relation (\ref{globalrelation}) to compute $J(t)$:
\begin{align}\nonumber
J(t) = & \int_{\partial D_3^0} \biggl\{-\frac{2}{\Delta}\biggl[kc - \Phi_1^{(1)} + \frac{\varphi_1^{(1)}\bar{\Phi}_2^{(1)}}{k}e^{2ikL} + \varphi_1^{(1)}e^{2ikL}\biggr] \biggr\}_- dk
	\\ \label{Jlaststep}
& + \int_{\partial D_3^0} \biggl\{\frac{2}{\Delta}\frac{\varphi_1^{(1)}\bar{\Phi}_2^{(1)}}{k}e^{2ikL} + \frac{2}{\Delta}[k\Phi_1 - \Phi_1^{(1)}]\biggr\}_-dk
	\\ \nonumber
& + \int_{\partial D_3^0} \frac{2k}{\Delta}\bigl[\Phi_1 (\bar{\varphi}_2 - 1) - (\bar{\Phi}_2 -1)\varphi_1e^{2ikL}\bigr]_-dk.
\end{align}
The asymptotics (\ref{casymptotics}) of $c(t,k)$ and Cauchy's theorem imply that the first integral on the rhs of (\ref{Jlaststep}) equals $- i\pi \varphi_1^{(2)}(t)$. 
Moreover, we can write the second integral on the rhs of (\ref{Jlaststep}) as
\begin{align*}
    \int_{\partial D_3^0} \biggl\{\frac{h_0}{2i}\Phi_{2-} + \frac{2}{\Delta}[k\Phi_{1-} + ig_0]\biggr\} dk.
\end{align*}
Therefore, equations (\ref{ipivarphi12}) and (\ref{Jlaststep}) imply
\begin{align}\nonumber
2 i\pi \varphi_1^{(2)}(t) = &\; \int_{\partial D_3^0} \biggl\{-\frac{k \Sigma \varphi_1}{\Delta} + \frac{\Sigma}{\Delta} \varphi_1^{(1)} \biggr\}_- dk
+  \int_{\partial D_3^0} \biggl\{\frac{h_0}{2i}\Phi_{2-} + \frac{2}{\Delta}[k\Phi_{1-} + ig_0]\biggr\} dk
	\\ \label{varphi12expression}
& + \int_{\partial D_3^0} \frac{2k}{\Delta}\bigl[\Phi_1 (\bar{\varphi}_2 - 1) - (\bar{\Phi}_2 -1)\varphi_1e^{2ikL}\bigr]_-dk.
\end{align}
Using (\ref{ipivarphi21}) and (\ref{varphi12expression}) in the second of equations (\ref{g1h1Phi}), we find (\ref{h1expression}).

$(b)$ In order to derive (\ref{g0expression}) we note that the first of equations (\ref{g0h0Phi}) expresses $g_0$ in terms of $\Phi_1^{(1)}$. Furthermore, equations (\ref{asymptoticformulas}) and Cauchy's theorem imply
\begin{align}
  -\frac{i\pi}{2} \Phi_1^{(1)}(t) = \int_{\partial D_2} \Phi_1 dk = \int_{\partial D_4} \Phi_1 dk.
\end{align}
Thus,
\begin{align}\label{ipiPhi11}
  i\pi \Phi_1^{(1)}(t) & = \biggl(\int_{\partial D_3} + \int_{\partial D_1}\biggr) \Phi_1 dk 
  = \int_{\partial D_3} \Phi_{1-} dk 
  = \int_{\partial D_3^0} \frac{\Sigma}{\Delta} \Phi_{1+}dk + I(t),
\end{align}
where the function $I(t)$ is defined by
$$I(t) = - \int_{\partial D_3^0} \frac{2}{\Delta} (e^{-2ikL} \Phi_1)_+ dk.$$
The last step involves using the global relation (\ref{globalrelation}) to compute $I(t)$:
\begin{align}\label{Neumannlaststep}
I(t) & =  - \int_{\partial D_3^0} \frac{2}{\Delta} (e^{-2ikL} c)_+ dk
+ \int_{\partial D_3^0} \frac{2}{\Delta} \bigl[\Phi_1(\bar{\varphi}_2 - 1)e^{-2ikL} - \bar{\Phi}_2\varphi_1\bigr]_+ dk .
\end{align}
The asymptotics (\ref{casymptotics}) of $c(t,k)$ and Cauchy's theorem imply that the first term on the rhs equals $-i\pi\Phi_1^{(1)}$; equations (\ref{ipiPhi11}), (\ref{Neumannlaststep}), and the first of equations (\ref{g0h0Phi}) yield (\ref{g0expression}). The proof of (\ref{h0expression}) is similar.
\proofend

\subsection{Effective characterizations}\label{effectivesubsec}
The substitution of the expressions (\ref{g1h1expressions}) for $g_1(t)$, $h_1(t)$ into (\ref{Phivarphieqs}) yields a system of quadratically nonlinear integral equations for $\{\Phi_j(t, k), \varphi_j(x,t)\}_1^2$. This nonlinear system provides an effective characterization of the spectral functions for the Dirichlet problem. In particular, given the Dirichlet data $g_0(t)$ and $h_0(t)$, the system can be solved recursively to all orders in a well-defined perturbative scheme. Indeed, substituting into (\ref{g1h1expressions}) the expansions
\begin{align} \nonumber
 &\Phi_j = \Phi_{j0} + \epsilon \Phi_{j1} + \epsilon^2 \Phi_{j2} + \cdots, &&
 \varphi_j = \varphi_{j0} + \epsilon \varphi_{j1} + \epsilon^2 \varphi_{j2} + \cdots, \quad j = 1,2,
  	\\ \nonumber
& g_0 = \epsilon g_{01} + \epsilon^2 g_{02} + \cdots, &&  g_1 = \epsilon g_{11} + \epsilon^2 g_{12} + \cdots, 
 	\\ \label{epsilonexpansions}
& h_0 = \epsilon h_{01} + \epsilon^2 h_{02} + \cdots, &&  h_1 = \epsilon h_{11} + \epsilon^2 h_{12} + \cdots, 
\end{align}
where $\epsilon > 0$ is a small parameter, the terms of $O(\epsilon^n)$ yield
\begin{subequations}\label{g1nh1n}
\begin{align}\label{g1nh1na}
g_{1n} = \frac{4}{\pi i}&\int_{\partial D_3^0} \biggl\{\frac{\Sigma}{2\Delta}\left[k\Phi_{1n-} + ig_{0n} \right] - \frac{1}{\Delta}\left[k\varphi_{1n-} + ih_{0n}\right] \biggr\}dk + \lot,
	\\ \label{g1nh1nb}
h_{1n} = \frac{4}{\pi i}&\int_{\partial D_3^0} \biggl\{ -\frac{\Sigma}{2\Delta} \left[k\varphi_{1n-} + ih_{0n}\right]
+ \frac{1}{\Delta}\left[k\Phi_{1n-} + ig_{0n}\right] \biggr\}dk + \lot.
\end{align}
\end{subequations}
The terms of $O(\epsilon^n)$ of the first and third equations in (\ref{Phivarphieqs}) yield
\begin{subequations}\label{Phi1nvarphi1n}
\begin{align}
&\Phi_{1n}(t, k) = \int_0^t e^{4ik^2(t' - t)}  (2kg_{0n}(t') + ig_{1n}(t')) dt' + \lot,
	\\
&  \varphi_{1n}(t, k) = \int_0^t e^{4ik^2(t' - t)} (2kh_{0n}(t') + ih_{1n}(t')) dt' + \lot.
\end{align}
\end{subequations}
The odd parts of the latter two equations yield
\begin{subequations}\label{Phivarphi1nminus}
\begin{align}
&\Phi_{1n-}(t, k) = 4k \int_0^t e^{4ik^2(t' - t)} g_{0n}(t')dt' + \lot,
	\\
&  \varphi_{1n-}(t, k) = 4k\int_0^t e^{4ik^2(t' - t)} h_{0n}(t') dt' + \lot.
\end{align}
\end{subequations}
It follows from (\ref{Phivarphi1nminus}) that $\Phi_{1n-}$ and $\varphi_{1n-}$ can be determined at each step from the known Dirichlet boundary values $g_{0n}$ and $h_{0n}$; $g_{1n}$ and $h_{1n}$ can then be determined from (\ref{g1nh1n}).

Similarly, the nonlinear system obtained by substituting the expressions (\ref{g0h0expressions}) for $g_0(t)$, $h_0(t)$ into (\ref{Phivarphieqs}) provides an effective characterization of the spectral functions for the Neumann problem.
Indeed, the terms of $O(\epsilon^n)$ of (\ref{g0expression})-(\ref{h0expression}) yield
\begin{subequations}\label{g0nh0n}
\begin{align}
 g_{0n}(t) =& \; \frac{2}{\pi} \int_{\partial D_3^0} \frac{1}{\Delta} \biggl\{ \frac{\Sigma}{2 }\Phi_{1n+} - \varphi_{1n+}\biggr\}dk,
	\\ 
h_{0n}(t) = & \; \frac{2}{\pi} \int_{\partial D_1^0} \frac{1}{\Delta}\biggl\{ - \frac{\Sigma}{2}\varphi_{1n+}   + \Phi_{1n+}\biggr\} dk,
\end{align}
\end{subequations}
while the even parts of (\ref{Phi1nvarphi1n}) yield
\begin{align*}
&\Phi_{1n+}(t, k) = 2i \int_0^t e^{4ik^2(t' - t)} g_{1n}(t') dt' + \lot,
	\\
&  \varphi_{1n+}(t, k) = 2i \int_0^t e^{4ik^2(t' - t)} h_{1n}(t') dt' + \lot.
\end{align*}
At each step in the perturbative scheme the functions $\Phi_{1n+}$ and $\varphi_{1n+}$ can be determined from the latter two equations and then $g_{0n}$ and $h_{0n}$ can be found from (\ref{g0nh0n}).

\subsection{The linear limit}
The linear limit of (\ref{g1expression}) yields $g_1 = \epsilon g_{11} + O(\epsilon^2)$ where
$$g_{11} =  
\frac{4}{\pi i}\int_{\partial D_3^0} \biggl\{\frac{\Sigma}{2\Delta}\left[k\Phi_{11-} + ig_{0n} \right] - \frac{1}{\Delta}\left[k\varphi_{11-} + ih_{0n}\right] \biggr\}dk.$$
Equation (\ref{Phivarphi1nminus}) becomes
$$\Phi_{11-} = 4k \int_0^t e^{4ik^2(s - t)} g_{01}(s) ds, \qquad
\varphi_{11-} = 4k \int_0^t e^{4ik^2(s - t)} h_{01}(s) ds.$$
Thus,
\begin{align}\label{g11formula}
g_{11} =  \frac{4}{\pi i}\int_{\partial D_3^0} \biggl\{\frac{\Sigma}{\Delta}\left(2k^2 \int_0^t e^{4ik^2(s - t)} g_{01} ds - \frac{g_{01}}{2i} \right) - \frac{4k^2}{\Delta} \int_0^t e^{4ik^2(s - t)} h_{01} ds + \frac{h_{01}}{i\Delta} \biggr\}dk.
\end{align}
This coincides with the formula of appendix \ref{linearlimitapp}, where the linearized equation $iq_t + q_{xx} = 0$ is solved directly.

\subsection{The large $L$ limit}
In the limit $L \to \infty$, the representations for $g_1$ and $g_0$ of theorem \ref{th1} reduce to the corresponding representations on the half-line. Indeed, as $L \to \infty$,
\begin{align*}
& h_0 \to 0, \qquad h_1 \to 0, \qquad \varphi_1 \to 0, \qquad \varphi_2 \to 1,
	\\
&\frac{\Sigma}{\Delta} \to 1  \text{ as $k \to \infty$ in $D_3$}, \qquad \frac{\Sigma}{\Delta} \to -1  \text{ as $k \to \infty$ in $D_1$}.
\end{align*}
%$$\mathcal{L}_j \to 0, \quad \mathcal{M}_j \to 0, \quad F(t,k) \to - \frac{i}{2}g_0(t)\hat{M}_2.$$
Thus, the $L \to \infty$ limits of the representations (\ref{g1expression}) and (\ref{g0expression}) are
\begin{align*}
g_1 =  \frac{4}{\pi i}\int_{\partial D_3} \biggl\{\frac{k}{2}\Phi_{1-} - \frac{g_0}{2i}  + \frac{ig_0}{2} \Phi_{2-} \biggr\}dk
 \quad \text{and}\quad
 g_0 = \frac{1}{\pi} \int_{\partial D_3^0} \Phi_{1+} dk,
\end{align*}
respectively, and these formulas coincide with the corresponding half-line formulas, cf. \cite{FLnonlinearizable}.

\section{The GLM approach}\nequation\label{glmsec}
In theorem \ref{th1} we derived effective representations for $\{g_1, h_1, g_0, h_0\}$ in terms of the eigenfunctions $\{\Phi_j, \varphi_j\}_1^2$. In what follows we will express the above boundary values in terms of the GLM representations.

For a function $f(t,s)$, we let $\hat{f}(t,k)$ denote the transform
$$\hat{f}(t,k) = \int_{-t}^t e^{2ik^2(s-t)} f(t,s) ds.$$
The eigenfunctions $\{\Phi_j, \varphi_j\}_1^2$ admit the following GLM representations
\begin{subequations}\label{GLMreps}
\begin{align}
& \Phi_1(t,k) = \hat{L}_1 - \frac{i}{2}g_0(t)\hat{M}_2 + k\hat{M}_1, &&
\Phi_2(t,k) = 1 + \hat{L}_2 + \frac{i\lambda}{2}\bar{g}_0 \hat{M}_1 + k\hat{M}_2,
	\\
& \varphi_1(t,k) = \hat{\mathcal{L}}_1 - \frac{i}{2}h_0(t)\hat{\mathcal{M}}_2 + k\hat{\mathcal{M}}_1, &&
\varphi_2(t,k) = 1 + \hat{\mathcal{L}}_2 + \frac{i\lambda}{2}\bar{h}_0 \hat{\mathcal{M}}_1 + k\hat{\mathcal{M}}_2,
\end{align}
\end{subequations}	
where the functions $\{L_j(t,s), M_j(t,s), \mathcal{L}_j(t,s), \mathcal{M}_j(t,s)\}_1^2$, $-t < s <t$, satisfy a nonlinear Goursat system (see \cite{FI2004}) together with the initial conditions
\begin{align}\nonumber
& L_1(t,t) = \frac{i}{2}g_1(t), \quad M_1(t,t) = g_0(t), \quad
\mathcal{L}_1(t,t) = \frac{i}{2}h_1(t), \quad \mathcal{M}_1(t,t) = h_0(t),
	\\ \label{LMinitialconditions}
& L_2(t,-t) = 	M_2(t,-t) = \mathcal{L}_2(t,-t) = \mathcal{M}_2(t,-t) = 0.
\end{align}

%\begin{align*}
%& \hat{L}_j(t,k) = \int_{-t}^t e^{2ik^2(s-t)} L_j(t,s) ds = 2\int_0^t e^{4ik^2(\tau -t)} L_j(t, 2\tau -t) d\tau,
%	\\
%& \hat{M}_j(t,k) = \int_{-t}^t e^{2ik^2(s-t)} M_j(t,s) ds = 2\int_0^t e^{4ik^2(\tau -t)} M_j(t, 2\tau -t) d\tau, \qquad j = 1,2.
%\end{align*}

\begin{theorem}\label{GLMth}
Define the function $F(t,k)$ by
\begin{align}\label{Fdef}
F(t,k) = \; & 
\frac{i}{2}h_0(t)e^{2ikL}\hat{\mathcal{M}}_2 - \frac{i}{2}g_0(t)\hat{M}_2 
 	\\ \nonumber
& +  \left(\overline{\hat{\mathcal{L}}}_2 - i\lambda \frac{h_0(t)}{2}\overline{\hat{\mathcal{M}}}_1  + k \overline{\hat{\mathcal{M}}}_2\right)\left(\hat{L}_1 - \frac{i}{2}g_0(t)\hat{M}_2 + k\hat{M}_1\right)
	\\\nonumber
& - e^{2ikL}\left(\overline{\hat{L}}_2 - i\lambda \frac{g_0(t)}{2}\overline{\hat{M}}_1  + k \overline{\hat{M}}_2\right)
 \left(\hat{\mathcal{L}}_1 - \frac{i}{2}h_0(t)\hat{\mathcal{M}}_2 + k\hat{\mathcal{M}}_1\right),
 \end{align}
where $\overline{\hat{\mathcal{L}}}_2$ is short-hand notation for $\overline{\hat{\mathcal{L}}_2(t, \bar{k})}$ etc.
Under the assumptions of theorem \ref{th1}, the following formulas are valid:
\begin{itemize}
\item[$(a)$] For the Dirichlet problem, the unknown boundary values $g_1$ and $h_1$ are given by
\begin{subequations}\label{GLMg1h1expressions}
\begin{align}\nonumber g_1(t) = \;&
\frac{4}{i\pi}  \int_{\partial D_1^0}\biggl\{ -\frac{2k^2}{\Delta} \left[\hat{\mathcal{M}}_1(t, k) - \frac{h_0(t)}{2ik^2}\right]
 +  \frac{k^2\Sigma}{\Delta} \left[\hat{M}_1(t, k) - \frac{g_0(t)}{2ik^2}\right]
	\\ \label{GLMg1expression}
& + \frac{k}{\Delta} (e^{-2ikL}F(t,k))_- 
 + \frac{ig_0(t)}{2} k \hat{M}_2(t, k) + \frac{k g_0(t)}{2i}\overline{\hat{\mathcal{M}}_2(t,\bar{k})}\biggr\} dk,
	\\ \nonumber
h_1(t) = \;&
 \frac{4}{i \pi} \int_{\partial D_1^0} \biggl\{\frac{2k^2}{\Delta} \left[\hat{M}_1(t,k) - \frac{g_0(t)}{2ik^2}\right]
 - \frac{k^2\Sigma}{\Delta} \left[\hat{\mathcal{M}}_1(t, k) - \frac{h_0(t)}{2ik^2}\right]
	\\ \label{GLMh1expression}
& + \frac{k}{\Delta}(F(t,k))_-
 + \frac{ih_0(t)}{2}  k \hat{\mathcal{M}}_2(t, k) +  \frac{k h_0(t)}{2i}\overline{\hat{M}_2(t,\bar{k})} \biggr\} dk.
\end{align}
\end{subequations}

\item[$(b)$] For the Neumann problem, the unknown boundary values $g_0$ and $h_0$ are given by 
\begin{subequations}\label{GLMg0h0expressions}
 \begin{align}\label{GLMg0expression}
& g_0(t) =\frac{2}{\pi} \int_{\partial D_1^0} \frac{1}{\Delta(k)} \left\{\Sigma(k)\hat{L}_1(t,k) - 2\hat{\mathcal{L}}_1(t,k)  + (e^{-2ikL}F(t,k))_+\right\}dk,
	\\ \label{GLMh0expression}
& h_0(t) =  \frac{2}{\pi} \int_{\partial D_1^0} \frac{1}{\Delta(k)}  \left\{- \Sigma(k)\hat{\mathcal{L}}_1(t,k)  + 2\hat{L}_1(t,k)  + F_+(t,k) \right\} dk.	
\end{align}
\end{subequations}
\end{itemize}
\end{theorem}

\begin{remark}\upshape
The representations in (\ref{GLMg1h1expressions}) coincide with the representations (4.7)-(4.8) in \cite{FI2004}, except that the last two terms on the rhs of (\ref{GLMg1expression}), as well as the last two terms on the rhs of (\ref{GLMh1expression}) were missed in  \cite{FI2004}. These terms, which arise from somewhat subtle boundary effects, are needed in order for equations (\ref{GLMg1h1expressions}) to be consistent with the representations of theorem \ref{th1} and are also required in order to obtain the correct large $L$ limit.
\end{remark}
  
\proofbegin
Let us first consider the Dirichlet problem. In view of the GLM representations (\ref{GLMreps}), we may write the global relation (\ref{globalrelation}a) as
\begin{align}\label{411FI2004}
-\hat{L}_1 + e^{2ikL}\hat{\mathcal{L}}_1
= k\hat{M}_1 - ke^{2ikL}\hat{\mathcal{M}}_1 + F  -c,
\end{align}
where
\begin{align}\label{intervalFexpression}
F(t,k) = 
\frac{i}{2}h_0(t)e^{2ikL}\hat{\mathcal{M}}_2 - \frac{i}{2}g_0(t)\hat{M}_2 
 +  (\bar{\varphi}_2 -1)\Phi_1 - e^{2ikL}(\bar{\Phi}_2 - 1)\varphi_1.
\end{align} 
The expression of $F$ can be expressed as in (\ref{Fdef}).
Letting $k \to -k$ in (\ref{411FI2004}), we find
\begin{align}\label{411FI2004b}
-\hat{L}_1 + e^{-2ikL}\hat{\mathcal{L}}_1
= -k\hat{M}_1 + ke^{-2ikL}\hat{\mathcal{M}}_1 + F(t,-k) - c(t,-k).
\end{align}
Solving (\ref{411FI2004}) and (\ref{411FI2004b}) for $\hat{L}_1$ and $\hat{\mathcal{L}}_1$, we find
\begin{align}
 \hat{\mathcal{L}}_1 & = \frac{2k}{\Delta}\hat{M}_1 - \frac{k\Sigma}{\Delta}\hat{\mathcal{M}}_1
  + \frac{1}{\Delta}(F - c)_-,
	\\
- \hat{L}_1 & = \frac{2k}{\Delta}\hat{\mathcal{M}}_1 - \frac{k\Sigma}{\Delta}\hat{M}_1
  - \frac{1}{\Delta}(e^{-2ikL}(F - c))_-.
\end{align}
Multiplying these equations by $k e^{4ik^2(t-t')}$, $0 < t' < t$, and integrating along $\partial D_1^0$ with respect to $dk$, we obtain
\begin{subequations}\label{glmderivationmidstep}
\begin{align} \nonumber
\int_{\partial D_1^0} k e^{4ik^2(t-t')}  \hat{\mathcal{L}}_1dk = \;&
\int_{\partial D_1^0} \frac{2k^2}{\Delta}e^{4ik^2(t-t')}\hat{M}_1 dk
- \int_{\partial D_1^0} \frac{k^2\Sigma}{\Delta} e^{4ik^2(t-t')}\hat{\mathcal{M}}_1 dk
	\\
& + \int_{\partial D_1^0} \frac{k}{\Delta}e^{4ik^2(t-t')}F_-dk,
	\\ \nonumber
- \int_{\partial D_1^0} k e^{4ik^2(t-t')} \hat{L}_1dk = \;&
\int_{\partial D_1^0}  \frac{2k^2}{\Delta}e^{4ik^2(t-t')}\hat{\mathcal{M}}_1 dk
-\int_{\partial D_1^0} \frac{k^2\Sigma}{\Delta}e^{4ik^2(t-t')}\hat{M}_1dk
	\\
& -  \int_{\partial D_1^0}  \frac{k}{\Delta} e^{4ik^2(t-t')}(e^{-2ikL}F)_- dk,
\end{align}
\end{subequations}
where we have used that the functions
$$\frac{k}{\Delta}(c(t,k))_- \qquad \text{and}\qquad \frac{k}{\Delta}(e^{-2ikL}c(t,k))_-$$
are bounded and analytic in $D_1^0$, so that their contributions vanish by Jordan's lemma.

The next step is to take the limit $t' \uparrow t$ in (\ref{glmderivationmidstep}). This can be achieved by using the identities %(eqs (4.3) and (4.4) in FI2004)
\begin{align}\label{identity1}
& \int_{\partial D_1} k e^{4ik^2(t-t')} \hat{f}(t,k) dk =  
%2 \int_{\partial D_1^0} k\int_0^t e^{4ik^2(\tau -t')} f(t, 2\tau -t) d\tau dk = 
%& \int_{\partial D_1} k \int_{-t}^t e^{2ik^2(s-t')} f(t,s) ds dk   = 
\begin{cases} \frac{\pi}{2}f(t, 2t' - t), &\quad 0 < t' < t, \\
\frac{\pi}{4} f(t, t), &\quad 0 < t' = t,
\end{cases}
\end{align}
and
\begin{align} \label{identity2}
& \int_{\partial D_1^0} \frac{k^2}{\Delta}  e^{4ik^2(t-t')} \hat{f}(t,k) dk
% = 2 \int_{\partial D_1^0} k^2\int_0^t e^{4ik^2(\tau -t')} L_j(t, 2\tau -t) d\tau dk
	\\ \nonumber
& \qquad = 2 \int_{\partial D_1^0} \frac{k^2}{\Delta} \left[\int_0^{t'} e^{4ik^2(\tau -t')} f(t, 2\tau -t) d\tau - \frac{f(t, 2t' - t)}{4ik^2}\right]dk, \qquad 0 < t' < t.
\end{align}	
The identity (\ref{identity2}) is also valid if $ \frac{k^2}{\Delta} $ is replaced by $k^2$ or $ \frac{k^2\Sigma}{\Delta}$.
Utilizing these identities in (\ref{glmderivationmidstep}), we find
\begin{align*}
\frac{\pi}{2}\mathcal{L}_1(t, 2t' - t) = \;&
4 \int_{\partial D_1^0} \frac{k^2}{\Delta} \left[\int_0^{t'} e^{4ik^2(\tau -t')} M_1(t, 2\tau -t) d\tau - \frac{M_1(t, 2t' - t)}{4ik^2}\right]dk
	\\
& - 2 \int_{\partial D_1^0} \frac{k^2\Sigma}{\Delta} \left[\int_0^{t'} e^{4ik^2(\tau -t')} \mathcal{M}_1(t, 2\tau -t) d\tau - \frac{\mathcal{M}_1(t, 2t' - t)}{4ik^2}\right]dk
	\\
& + \int_{\partial D_1^0} \frac{k}{\Delta}e^{4ik^2(t-t')}F_-(t,k)dk
\end{align*}
and
\begin{align*}
- \frac{\pi}{2}L_1(t, 2t' - t) = \;&
4 \int_{\partial D_1^0} \frac{k^2}{\Delta} \left[\int_0^{t'} e^{4ik^2(\tau -t')} \mathcal{M}_1(t, 2\tau -t) d\tau - \frac{\mathcal{M}_1(t, 2t' - t)}{4ik^2}\right]dk
	\\
& -2 \int_{\partial D_1^0} \frac{k^2\Sigma}{\Delta} \left[\int_0^{t'} e^{4ik^2(\tau -t')} M_1(t, 2\tau -t) d\tau - \frac{M_1(t, 2t' - t)}{4ik^2}\right]dk
	\\
& -  \int_{\partial D_1^0}  \frac{k}{\Delta} e^{4ik^2(t-t')}(e^{-2ikL}F(t,k))_- dk.
\end{align*}
Letting $t' \uparrow t$ in these equations and using the initial conditions (\ref{LMinitialconditions}) as well as the following lemma, we find the representations in (\ref{GLMg1h1expressions}).

\begin{lemma}\label{Fintegrallimitlemma}
\begin{subequations}\label{Fintegrallimit}
\begin{align} \label{Fintegrallimita}
& \lim_{t' \uparrow t} \int_{\partial D_1^0} \frac{k}{\Delta}e^{4ik^2(t-t')}F_-(t,k)dk
 =  \int_{\partial D_1^0} \frac{k}{\Delta}F_-(t,k)dk
	\\ \nonumber
& \hspace{3cm} + \frac{ih_0(t)}{2} \int_{\partial D_1^0} k \hat{\mathcal{M}}_2(t, k)dk
+ \int_{\partial D_1^0} \frac{k h_0(t)}{2i}\overline{\hat{M}_2(t, \bar{k})} dk,	
 	\\ \label{Fintegrallimitb}
& \lim_{t' \uparrow t} \int_{\partial D_1^0} \frac{k}{\Delta} e^{4ik^2(t-t')}(e^{-2ikL}F(t,k))_- dk
=  \int_{\partial D_1^0} \frac{k}{\Delta} (e^{-2ikL}F(t,k))_- dk
	\\ \nonumber
& \hspace{3cm} + \frac{ig_0(t)}{2} \int_{\partial D_1^0} k \hat{M}_2(t, k)dk
+ \int_{\partial D_1^0} \frac{k g_0(t)}{2i}\overline{\hat{\mathcal{M}}_2(t,\bar{k})} dk.
\end{align}
\end{subequations}
\end{lemma}
\proofbegin
We prove (\ref{Fintegrallimitb}); the proof of (\ref{Fintegrallimita}) is similar. 
If one naively takes the limit inside the integral in (\ref{Fintegrallimitb}), one finds the first term on the rhs of (\ref{Fintegrallimitb}). 
The other two terms on the rhs of (\ref{Fintegrallimitb}) arise from interchanging the limit and the integration. We will next describe how these terms arise in detail.

We write
\begin{align} \nonumber
\int_{\partial D_1^0}& \frac{k}{\Delta} e^{4ik^2(t-t')}(e^{-2ikL}F(t,k))_- dk
=  \int_{\partial D_1^0} k e^{4ik^2(t-t')} \frac{ig_0}{2}\hat{M}_2 dk
	\\ \nonumber
& - \int_{\partial D_1^0} k e^{4ik^2(t-t')} \biggl(\overline{\hat{\mathcal{L}}}_2 - i\lambda \frac{h_0}{2} \overline{\hat{\mathcal{M}}}_1\biggr)\biggl(\hat{L}_1 - \frac{i}{2}g_0 \hat{M}_2\biggr) dk
	\\ \nonumber
& - 2 \int_{\partial D_1^0} \frac{k^2}{\Delta} e^{4ik^2(t-t')} \left[\biggl(\overline{\hat{L}}_2 - \frac{i\lambda g_0}{2}\overline{\hat{M}}_1\biggr)\hat{\mathcal{M}}_1
+ \overline{\hat{M}}_2\biggl(\hat{\mathcal{L}}_1 - \frac{ih_0}{2} \hat{\mathcal{M}}_2\biggr)\right] dk
	\\ \nonumber
& +\int_{\partial D_1^0} \frac{k^2\Sigma}{\Delta} e^{4ik^2(t-t')}
\biggl[\biggl(\overline{\hat{\mathcal{L}}}_2 - \frac{i\lambda h_0}{2}\overline{\hat{\mathcal{M}}}_1\biggr) \hat{M}_1
 + \biggl(\hat{L}_1 - \frac{ig_0}{2}\hat{M}_2\biggr)\overline{\mathcal{\hat{M}}}_2\biggr]dk
	\\  \label{Fintegralexpanded}
& - \int_{\partial D_1^0} k^3 \overline{\hat{\mathcal{M}}}_2\hat{M}_1 e^{4ik^2(t-t')}dk.
\end{align}
The first integral on the rhs of (\ref{Fintegralexpanded}) yields the following contribution in the limit $t' \to t$:
$$\lim_{t' \uparrow t} \frac{ig_0(t)}{2} \int_{\partial D_1^0} k e^{4ik^2(t-t')} \hat{M}_2(t,k) dk
=\lim_{t' \uparrow t}  \frac{ig_0(t)}{2} \frac{\pi}{2} M_2(t, 2t' - t) = \frac{i\pi g_0(t)}{4}M_2(t,t).$$
On the other hand, utilizing the second row of (\ref{identity1}),
$$\frac{ig_0(t)}{2} \int_{\partial D_1^0} k \hat{M}_2(t,k) dk
= \frac{i \pi g_0(t)}{8}M_2(t,t).$$
Therefore,
\begin{align}\label{limtptot}
\lim_{t' \uparrow t} \frac{ig_0(t)}{2} \int_{\partial D_1^0} k e^{4ik^2(t-t')} \hat{M}_2(t,k) dk
=\;& \frac{ig_0(t)}{2} \int_{\partial D_1^0} k \hat{M}_2(t,k) dk
	\\ \nonumber
& + \frac{ig_0(t)}{2} \int_{\partial D_1^0} k \hat{M}_2(t, k)dk.
\end{align}
The first term on the rhs of (\ref{limtptot}) is the contribution obtained by taking the limit inside the integral; this term is included in the first term on the rhs of (\ref{Fintegrallimitb}). In addition to this term, there is also an additional term arising from the interchange of the limit and the integration; this is the second term on the rhs of (\ref{Fintegrallimitb}). 

We now consider the last integral on the rhs of (\ref{Fintegralexpanded}), which can be written as
\begin{align}\label{thirdintegralinclaim}
&- \int_{\partial D_1^0} k^3 \overline{\hat{\mathcal{M}}}_2\hat{M}_1 e^{4ik^2(t-t')}dk
= -2\int_{\partial D_1^0} k^3 \overline{\hat{\mathcal{M}}_2(t,\bar{k})}\int_0^t e^{4ik^2(\tau -t')} M_1(t, 2\tau -t) d\tau dk.
\end{align}
The rhs of (\ref{thirdintegralinclaim}) equals
\begin{align}\label{thirdintegralrhs}
-2 \int_{\partial D_1^0} k^3 \overline{\hat{\mathcal{M}}_2(t,\bar{k})}
\left\{\int_0^{t'} e^{4ik^2(\tau -t')} M_1(t, 2\tau -t) d\tau  -  \frac{M_1(t, 2t' - t)}{4ik^2}\right\}dk.
\end{align}
Indeed, the rhs of (\ref{thirdintegralinclaim}) equals this term plus the following expression:
$$-2 \int_{\partial D_1^0} k^3 \overline{\hat{\mathcal{M}}_2(t,\bar{k})}
\left\{\int_{t'}^{t} e^{4ik^2(\tau -t')} M_1(t, 2\tau -t) d\tau  +  \frac{M_1(t, 2t' - t)}{4ik^2}\right\}dk.$$
Integration by parts shows that this integral vanishes by Jordan's lemma, because $\hat{\mathcal{M}}_2(t,k)$ is of $O(1/k^2)$ as $k \to \infty$ within $D_4$. 
Taking the limit $t' \uparrow t$ in (\ref{thirdintegralrhs}) and using (\ref{LMinitialconditions}), we find that the contribution of the last integral in (\ref{Fintegralexpanded}) is
$$- \lim_{t' \uparrow t}\int_{\partial D_1^0} k^3 \overline{\hat{\mathcal{M}}}_2\hat{M}_1 e^{4ik^2(t-t')}dk
=- \int_{\partial D_1^0} k^3 \overline{\hat{\mathcal{M}}}_2\hat{M}_1dk
+ \int_{\partial D_1^0} \frac{k g_0(t)}{2i}\overline{\hat{\mathcal{M}}}_2 dk.$$
The first term on the rhs is the contribution obtained by taking the limit inside the integral.
In addition to this term, there is also an additional term arising from the interchange of the limit and the integration; this is the third term on the rhs of (\ref{Fintegrallimitb}). 

Finally, we claim that the limits of the second, third, and fourth integrals on the rhs of (\ref{Fintegralexpanded}) can be computed by simply taking the limit inside the integral, i.e. in these cases no additional terms arise. We show this for the term 
$$I := \int_{\partial D_1^0} \frac{k^2\Sigma}{\Delta} e^{4ik^2(t-t')}
\biggl(\overline{\hat{\mathcal{L}}}_2 - \frac{i\lambda h_0}{2}\overline{\hat{\mathcal{M}}}_1\biggr) \hat{M}_1dk;
$$
the proofs for the other terms are similar.
We have
$$I = 2\int_{\partial D_1^0} \frac{k^2\Sigma}{\Delta}
\biggl(\overline{\hat{\mathcal{L}}}_2 - \frac{i\lambda h_0}{2}\overline{\hat{\mathcal{M}}}_1\biggr)
\int_0^t e^{4ik^2(\tau - t')}M_1(t, 2\tau -t)d\tau dk.$$
We can write this as
$$I = 2\int_{\partial D_1^0} \frac{k^2\Sigma}{\Delta}
\biggl(\overline{\hat{\mathcal{L}}}_2 - \frac{i\lambda h_0}{2}\overline{\hat{\mathcal{M}}}_1\biggr)
\left[\int_0^{t'} e^{4ik^2(\tau - t')}M_1(t, 2\tau -t)d\tau - \frac{M_1(t, 2t' -t)}{4ik^2}\right]dk.$$
Indeed, the difference between the preceding two expressions, 
$$2\int_{\partial D_1^0} \frac{k^2\Sigma}{\Delta}
\biggl(\overline{\hat{\mathcal{L}}}_2 - \frac{i\lambda h_0}{2}\overline{\hat{\mathcal{M}}}_1\biggr)
\left[\int_{t'}^t e^{4ik^2(\tau - t')}M_1(t, 2\tau -t)d\tau + \frac{M_1(t, 2t' -t)}{4ik^2}\right]dk,$$
vanishes by integration by parts and Jordan's lemma, since $\hat{\mathcal{L}}_2$ and $\hat{\mathcal{M}}_1$ are of $O(1/k^2)$ as $k \to \infty$ in $D_4$.
Taking the limit $t' \uparrow t$, we find
\begin{align*}
&\lim_{t' \uparrow t}
\int_{\partial D_1^0} \frac{k^2\Sigma}{\Delta} e^{4ik^2(t-t')}
\biggl(\overline{\hat{\mathcal{L}}}_2 - \frac{i\lambda h_0}{2}\overline{\hat{\mathcal{M}}}_1\biggr) \hat{M}_1dk
	\\
&= \int_{\partial D_1^0} \frac{k^2\Sigma}{\Delta}
\biggl(\overline{\hat{\mathcal{L}}}_2 - \frac{i\lambda h_0}{2}\overline{\hat{\mathcal{M}}}_1\biggr)
\left[\hat{M}_1(t, k) - \frac{g_0(t)}{2ik^2}\right]dk.
\end{align*}
However, in this case the additional term
$$-\int_{\partial D_1^0} \frac{k^2\Sigma}{\Delta}
\biggl(\overline{\hat{\mathcal{L}}}_2 - \frac{i\lambda h_0}{2}\overline{\hat{\mathcal{M}}}_1\biggr)\frac{g_0(t)}{2ik^2}dk,$$
vanishes because the integrand is analytic and of $O(1/k^2)$ as $k \to \infty$ in $D_1$. This completes the proof of lemma \ref{Fintegrallimitlemma}.
\proofend

We now return to theorem \ref{GLMth} and consider the Neumann problem. Solving (\ref{411FI2004}) and (\ref{411FI2004b}) for $\hat{M}_1$ and $\hat{\mathcal{M}}_1$, we find
\begin{align}
&k\hat{\mathcal{M}}_1 = \frac{2}{\Delta}\hat{L}_1 - \frac{\Sigma}{\Delta}\hat{\mathcal{L}}_1 + \frac{1}{\Delta}(F - c)_+, 
	\\
&k\hat{M}_1 = \frac{\Sigma}{\Delta}\hat{L}_1 - \frac{2}{\Delta}\hat{\mathcal{L}}_1 + \frac{1}{\Delta}(e^{-2ikL}(F - c))_+.\end{align}
Multiplying these equations by $k e^{4ik^2(t-t')}$, $0 < t' < t$, and integrating along $\partial D_1^0$ with respect to $dk$, 
we find
\begin{align*}
\frac{\pi}{2} \mathcal{M}_1(t, 2t' - t) = & \int_{\partial D_1^0}  e^{4ik^2(t-t')}\biggl\{
\frac{2}{\Delta}\hat{L}_1  - \frac{\Sigma}{\Delta}\hat{\mathcal{L}}_1 + \frac{1}{\Delta}F_+ \biggr\}dk, 
	\\
 \frac{\pi}{2}M_1(t, 2t' - t) = & \int_{\partial D_1^0}  e^{4ik^2(t-t')} \biggl\{ \frac{\Sigma}{\Delta}\hat{L}_1 dk
- \frac{2}{\Delta}\hat{\mathcal{L}}_1 + \frac{1}{\Delta}(e^{-2ikL}F)_+\biggr\} dk,
\end{align*}
where we used that the functions
$$\frac{1}{\Delta}(c(t,k))_+, \qquad \frac{1}{\Delta}(e^{-2ikL}c(t,k))_+$$
are bounded and analytic in the interior of $\partial D_1^0$ so that their contributions vanish by Jordan's lemma.
Letting $t' \uparrow t$ in these equations and using the initial conditions (\ref{LMinitialconditions}) as well as the following lemma, we find the representations in (\ref{GLMg0h0expressions}).

\begin{lemma}\label{lemma44}
\begin{subequations}
\begin{align} \label{Fplusintegrallimita}
& \lim_{t' \uparrow t} \int_{\partial D_1^0} \frac{1}{\Delta}e^{4ik^2(t-t')}F_+(t,k)dk
 = \int_{\partial D_1^0} \frac{1}{\Delta}F_+(t,k)dk,
 	\\ \label{Fplusintegrallimitb}
& \lim_{t' \uparrow t} \int_{\partial D_1^0} \frac{1}{\Delta} e^{4ik^2(t-t')}(e^{-2ikL}F(t,k))_+ dk
= \int_{\partial D_1^0} \frac{1}{\Delta} (e^{-2ikL}F(t,k))_+ dk.
\end{align}
\end{subequations}
\end{lemma}
\proofbegin
Note that
\begin{align*} \nonumber
&\int_{\partial D_1^0} \frac{1}{\Delta} e^{4ik^2(t-t')} F_+(t,k) dk
= \int_{\partial D_1^0} \frac{1}{\Delta} e^{4ik^2(t-t')} \biggl\{
\frac{i}{2}h_0(t)\Sigma \hat{\mathcal{M}}_2 - ig_0(t)\hat{M}_2 
 	\\
& +  2\left(\overline{\hat{\mathcal{L}}}_2 - i\lambda \frac{h_0}{2}\overline{\hat{\mathcal{M}}}_1\right)\left(\hat{L}_1 - \frac{i}{2}g_0(t)\hat{M}_2\right) +  2k^2 \overline{\hat{\mathcal{M}}}_2 \hat{M}_1
	\\
& - \Sigma \left(\overline{\hat{L}}_2 - i\lambda \frac{g_0}{2}\overline{\hat{M}}_1\right)
 \left(\hat{\mathcal{L}}_1 - \frac{i}{2}h_0(t)\hat{\mathcal{M}}_2\right)
	\\
& - \Delta k \overline{\hat{M}}_2 \left(\hat{\mathcal{L}}_1 - \frac{i}{2}h_0(t)\hat{\mathcal{M}}_2 \right)
 - \Delta  \left(\overline{\hat{L}}_2 - i\lambda \frac{g_0}{2}\overline{\hat{M}}_1\right) k\hat{\mathcal{M}}_1
 - \Sigma k^2\overline{\hat{M}}_2\hat{\mathcal{M}}_1\biggr\}.
\end{align*}
Equation (\ref{Fplusintegrallimita}) now follows by using arguments similar to those that led to lemma \ref{Fintegrallimitlemma}. The proof of (\ref{Fplusintegrallimitb}) is similar. This completes the proof of lemma \ref{lemma44} and hence of theorem \ref{GLMth}.
\proofend

\subsection{Equivalence of the two representations}
We will show that the representations derived using the GLM approach in theorem \ref{GLMth} coincide with those of theorem \ref{th1}. 

\subsubsection{The representations for $g_1$ and $h_1$}
Using the expression (\ref{intervalFexpression}) for $F$ as well as the formulas
$$\hat{M}_j = \frac{1}{2k}\Phi_{j-}, \qquad \hat{\mathcal{M}}_j = \frac{1}{2k}\varphi_{j-}, \qquad j = 1,2,$$
we can write the representation (\ref{GLMg1expression}) of $g_1$ as
\begin{align*}
 g_1(t) = 
\frac{4}{i\pi} \int_{\partial D_1^0} \biggl\{&-\frac{1}{\Delta} \left[k\varphi_{1-} + ih_0(t)\right]
+ \frac{\Sigma}{2\Delta} \left[k\Phi_{1-}(k) + ig_0(t)\right]
	\\ \nonumber
& + \frac{k}{\Delta} \left[(\bar{\varphi}_2 -1)\Phi_1e^{-2ikL} - (\bar{\Phi}_2 - 1)\varphi_1\right]_- 
	\\ \nonumber
&- \frac{k}{\Delta} \left[e^{-2ikL}\frac{ig_0}{4k}\Phi_{2-}\right]_- 
+ \frac{ig_0(t)}{4}\Phi_{2-} + \frac{g_0(t)}{4i}\bar{\varphi}_{2-} \biggr\} dk,
\end{align*}
The identities
$$\int_{\partial D_1^0} \overline{\varphi_{2-}(t,\bar{k})} dk 
= -\pi i\overline{\varphi_2^{(1)}} = \pi i\varphi_2^{(1)} 
= \int_{\partial D_1^0} \varphi_{2-} dk,$$
imply that the function $\bar{\varphi}_{2-}$ in the above integrand can be replaced with $\varphi_{2-}$.
Moreover, since the integrand is an odd function of $k$, the contour $\partial D_1^0$ can be replaced with $\partial D_3^0$. In view of the identity 
$$- \frac{k}{\Delta} \left[e^{-2ikL}\frac{ig_0}{4k}\Phi_{2-}\right]_- 
+ \frac{ig_0(t)}{4}\Phi_{2-} = \frac{ig_0(t)}{2}\Phi_{2-},$$
we find the representation for $g_1$ in (\ref{g1expression}).
Similar computations show that the representations for $h_1$ are also equivalent.

\subsubsection{The representations for $g_0$ and $h_0$}
Using the expression (\ref{intervalFexpression}) for $F$ as well as the formulas
$$\Phi_{1+} = 2\hat{L}_1 - ig_0 \hat{M}_2, \qquad 
%\Phi_{2+} =2+ 2\hat{L}_2 + i\lambda\bar{g}_0 \hat{M}_1,
\varphi_{1+} = 2\hat{\mathcal{L}}_1 - ih_0 \hat{\mathcal{M}}_2, 
%\qquad \varphi_{2+} =2+ 2\hat{\mathcal{L}}_2 + i\lambda\bar{h}_0 \hat{\mathcal{M}}_1,
$$
a straightforward computation shows that the representations (\ref{g0h0expressions}) and (\ref{GLMg0h0expressions}) are equivalent.

\appendix
\section{The linear limit} \label{linearlimitapp}
\renewcommand{\theequation}{A.\arabic{equation}}\nequation
In this appendix we analyze the linearized version of the NLS equation, $iq_t + q_{xx} = 0$, on the interval $[0,L]$ using the unified method of \cite{F1997}. The global relation for this equation is (cf. equation (2.8) in \cite{Fbook})
\begin{align}\label{linearGR}
  \hat{q}_0(k) - \tilde{g}(ik^2) + e^{-ikL}\tilde{h}(ik^2) = e^{ik^2T}\hat{q}_T(k), \qquad k \in \C,
\end{align}  
where
\begin{align*}
  & \hat{q}_0(k) = \int_0^L e^{-ikx} q_0(x) dx, &&
   \hat{q}_T(k) = \int_0^L e^{-ikx} q(x,T) dx, %\qquad w(k) = ik^2, \quad c_0(k) = -k, \quad c_1(k) = i,
  	\\
  & \tilde{g}(k) = -k\tilde{g}_0(ik^2) + i\tilde{g}_1(ik^2), &&
  \tilde{h}(k) = -k\tilde{h}_0(ik^2) + i\tilde{h}_1(ik^2), 
  	\\
  & \tilde{g}_j(ik^2) = \int_0^T e^{ik^2 s}g_j(s) ds, &&
  \tilde{h}_j(ik^2) = \int_0^T e^{ik^2 s}h_j(s) ds, \qquad j = 0,1,	
\end{align*}
and $g_j(t), h_j(t)$, $j = 0,1$, denote the Dirichlet and Neumann boundary values as in (\ref{g0h0g1h1def}).
Equation (\ref{linearGR}) and the equation obtained by letting $k \to -k$ in (\ref{linearGR}) are the following equations:
\begin{align*}
& \hat{q}_0(k) + k\tilde{g}_0(ik^2) - i\tilde{g}_1(ik^2) - k\tilde{h}_0(ik^2)e^{-ikL} + i\tilde{h}_1(ik^2)e^{-ikL} = e^{ik^2 T} \hat{q}_T(k),
	\\
& \hat{q}_0(-k) - k\tilde{g}_0(ik^2) - i\tilde{g}_1(ik^2) + k\tilde{h}_0(ik^2)e^{ikL} + i\tilde{h}_1(ik^2)e^{ikL} = e^{ik^2 T} \hat{q}_T(-k).
\end{align*}
Eliminating $\tilde{h}_1$ from these equations and then solving for $\tilde{g}_1$, we find
$$\tilde{g}_1(ik^2) = \frac{i}{\Delta(k/2)}\Bigl[e^{ik^2T}(e^{ikL} \hat{q}_T(k))_- - (e^{ikL}\hat{q}_0(k))_- - k\Sigma(k/2) \tilde{g}_0(ik^2) + 2k\tilde{h}_0(ik^2)\Bigr].$$
Multiplying this equation by $k e^{-ik^2t}$ and integrating along $\partial D_1^0$ with respect to $dk$, the term involving $\hat{q}_T$ is eliminated and we find
$$\pi g_1 = \int_{\partial D_1^0}\frac{i k e^{-ik^2t}}{\Delta(k/2)}\Bigl[- (e^{ikL}\hat{q}_0(k))_- - k\Sigma(k/2) \tilde{g}_0(ik^2) + 2k\tilde{h}_0(ik^2)\Bigr].$$
Assuming that $q_0 = 0$ 
%$$g_1 = \frac{1}{\pi i}\int_{\partial D_1} \left[ \frac{k^2\Sigma(k/2)}{\Delta(k/2)}  \int_0^T e^{ik^2(s -t)}g_0(s) ds - \frac{2k^2}{\Delta(k/2)} \int_0^T e^{ik^2(s-t)}h_0(s) ds\right]dk.$$
and performing the change of variables $k = -2l$, we find
$$g_1(t) = \frac{4}{\pi i}\int_{\partial D_3^0} \left[\frac{2l^2\Sigma(l)}{\Delta(l)}  \int_0^T e^{4il^2(s -t)}g_0(s) ds - \frac{4l^2}{\Delta(l)} \int_0^T e^{4il^2(s-t)}h_0(s) ds\right]dl.$$
Using the identity 
%(cf. equation (4.4) in \cite{FI2004})
\begin{align}\label{identity44}
\int_{\partial D_3^0} \frac{k^2}{\Delta(k)}& \int_0^T e^{4ik^2(s - t)} K(s,T) ds dk
	\\ \nonumber
& = \int_{\partial D_3^0} \frac{k^2}{\Delta(k)} \left[\int_0^t e^{4ik^2(s-t)}K(s, T) ds - \frac{K(t, T)}{4ik^2}\right]dk, \qquad 0 < t < T,
\end{align}
and a similar identity obtained by replacing $\frac{k^2}{\Delta(k)}$ by $\frac{k^2\Sigma(k)}{\Delta(k)}$ in (\ref{identity44}), we arrive at
\begin{align*}
&g_1(t) = \frac{4}{\pi i}\int_{\partial D_3^0} \biggl\{\frac{2l^2\Sigma(l)}{\Delta(l)}  \int_0^t e^{4il^2(s -t)}g_0(s) ds
- \frac{\Sigma(l) g_0}{2i\Delta(l) }
 	\\
& \hspace{3cm} - \frac{4l^2}{\Delta(l)} \int_0^t e^{4il^2(s-t)}h_0(s) ds + \frac{h_0}{i\Delta(l)}\biggr\}dl,
\end{align*}
which coincides with formula (\ref{g11formula}) for $g_{11}$.

\section{The asymptotics of $c(t,k)$} \label{clemmaapp}
\renewcommand{\theequation}{B.\arabic{equation}}\nequation
We will prove lemma \ref{clemma}. We will make no assumption on the function $c(t,k)$ except that it satisfies the global relation (\ref{GR}) and that it has the boundedness properties stated in (\ref{cO}). 
The functions $\{\Phi_j, \varphi_j\}_1^2$ satisfy the systems
\begin{align}\label{Phi12system}
 \begin{cases}
  \Phi_{1t} = -4ik^2\Phi_1 - i\lambda |g_0|^2 \Phi_1 + (2kg_0 + ig_1)\Phi_2,
  	\\ 
\Phi_{2t} = \lambda(2k\bar{g}_0 - i\bar{g}_1) \Phi_1 + \lambda i |g_0|^2 \Phi_2,
 \end{cases}
\end{align}
and
\begin{align}\label{varphi12system}
 \begin{cases}
  \varphi_{1t} = -4ik^2\varphi_1 - i\lambda |h_0|^2 \varphi_1 + (2kh_0 + ih_1)\varphi_2,
  	\\ 
\varphi_{2t} = \lambda(2k\bar{h}_0 - i\bar{h}_1) \varphi_1 + \lambda i |h_0|^2 \varphi_2.
 \end{cases}
\end{align}
In view of the initial conditions
$$\Phi_1(0,k) = \varphi_1(0,k) = 0, \qquad \Phi_2(0,k) = \varphi_2(0,k) = 1,$$
this leads to the asymptotic expansions
\begin{align*}
&\begin{pmatrix} \Phi_1(t,k) \\ \Phi_2(t,k) \end{pmatrix}
= \begin{pmatrix} 0 \\ 1 \end{pmatrix} + \begin{pmatrix} \Phi_1^{(1)}(t) \\ \Phi_2^{(1)}(t) \end{pmatrix} \frac{1}{k}
+ \begin{pmatrix} \Phi_1^{(2)}(t) \\ \Phi_2^{(2)}(t) \end{pmatrix} \frac{1}{k^2} + O\Bigl(\frac{1}{k^3}\Bigr)
	\\ \nonumber
&\qquad + \left[\begin{pmatrix} -\Phi_1^{(1)}(0) \\ 0 \end{pmatrix}\frac{1}{k} 
+  \begin{pmatrix} - \Phi_1^{(2)}(0) + \Phi_1^{(1)}(0)\int_{(0,0)}^{(0,t)} \omega \\ - \frac{i \lambda}{2}\bar{g}_0(t)\Phi_1^{(1)}(0) \end{pmatrix}\frac{1}{k^2} + O\Bigl(\frac{1}{k^3}\Bigr)\right]e^{-4ik^2t},
	\\ \nonumber
& \hspace{11cm} k \to \infty, \; k \in \C,
	\\ \nonumber
&\begin{pmatrix} \varphi_1(t,k) \\ \varphi_2(t,k) \end{pmatrix}
= \begin{pmatrix} 0 \\ 1 \end{pmatrix} + \begin{pmatrix} \varphi_1^{(1)}(t) \\ \varphi_2^{(1)}(t) \end{pmatrix} \frac{1}{k}
+ \begin{pmatrix} \varphi_1^{(2)}(t) \\ \varphi_2^{(2)}(t) \end{pmatrix} \frac{1}{k^2} + O\Bigl(\frac{1}{k^3}\Bigr)
	\\ \nonumber
&\qquad + \left[\begin{pmatrix} -\varphi_1^{(1)}(0) \\ 0 \end{pmatrix}\frac{1}{k} 
+  \begin{pmatrix} - \varphi_1^{(2)}(0) + \varphi_1^{(1)}(0)\int_{(L,0)}^{(L,t)} \omega \\ - \frac{i \lambda}{2}\bar{h}_0(t)\varphi_1^{(1)}(0) \end{pmatrix}\frac{1}{k^2} + O\Bigl(\frac{1}{k^3}\Bigr)\right]e^{-4ik^2t},
	\\ \nonumber
& \hspace{11cm} k \to \infty, \; k \in \C.
\end{align*}
Similarly, we have
\begin{align}
& \begin{pmatrix} b(k) \\ a(k) \end{pmatrix} = \begin{pmatrix} 0 \\ 1 \end{pmatrix} + \begin{pmatrix} b^{(1)} \\ a^{(1)} \end{pmatrix} \frac{1}{k}
+ \begin{pmatrix} b^{(2)} \\ a^{(2)} \end{pmatrix} \frac{1}{k^2} + O\Bigl(\frac{1}{k^3}\Bigr)
	\\ \nonumber
& \qquad+ \left[\begin{pmatrix} -\frac{q_0(L)}{2i} \\ 0 \end{pmatrix}\frac{1}{k} 
+  \begin{pmatrix} \beta^{(2)} \\ \alpha^{(2)} \end{pmatrix}\frac{1}{k^2} + O\Bigl(\frac{1}{k^3}\Bigr)\right]e^{2ikL}, \qquad k \to \infty, \; k \in \C.
\end{align}
%\beta^{(2)} = - \frac{q_x(L,0)}{4} +  \frac{q_0(L)}{2i}\int_{(L,0)}^{(0,0)} \omega
Substituting these expansions into the global relation (\ref{GR}), we find 
\begin{align}\nonumber
  c(t,k) = &\; \biggl\{O\Bigl(\frac{1}{k^2}\Bigr) + O\Bigl(\frac{1}{k^2}\Bigr)e^{2ikL}\biggr\}e^{4ik^2t} 
  	\\ \nonumber
& + \biggl\{\frac{b^{(1)} - \Phi_1^{(1)}(0)}{k} + O\Bigl(\frac{1}{k^2}\Bigr)   + \biggl(\frac{\varphi_1^{(1)}(0)-\frac{q_0(L)}{2i}}{k} + O\Bigl(\frac{1}{k^2}\Bigr)\biggr)e^{2ikL}
    \biggr\} e^{-4ik^2t}
  	\\\label{calmostthere}
&   + \frac{\Phi_1^{(1)}(t)}{k} + \frac{\Phi_1^{(2)}(t) + \Phi_1^{(1)}(t)\bigl(a^{(1)} + \bar{\varphi}_2^{(1)}(t)\bigr)}{k^2} + O \Bigl(\frac{1}{k^3} \Bigr)
 	\\ \nonumber
&- \biggl\{\frac{\varphi_1^{(1)}(t)}{k} + \frac{\varphi_1^{(2)}(t) + \varphi_1^{(1)}(t)\bigl(\bar{a}^{(1)} + \bar{\Phi}_2^{(1)}(t)\bigr)}{k^2} + O \Bigl(\frac{1}{k^3} \Bigr)\biggr\}e^{2ikL}, 
	\\ \nonumber
& \hspace{7cm} k \to \infty, \; k \in \C.
\end{align}
The assumption that $c(t,k)$ is of $O((1 + e^{2ikL})/k)$ as $k \to \infty$ implies that the terms in (\ref{calmostthere}) involving $e^{-4ik^2t}$ and $e^{4ik^2t}$ must vanish, i.e. for consistency we require
\begin{align}\label{b1Phi11}
b^{(1)} = \Phi_1^{(1)}(0),  \qquad
\frac{q_0(L)}{2i} = \varphi_1^{(1)}(0). 
\end{align}
Using the expressions for $b^{(1)}$, $ \Phi_1^{(1)}(0)$, and $\varphi_1^{(1)}(0)$, we find that equations (\ref{b1Phi11}) are valid iff the initial and boundary conditions are compatible, i.e. iff
\begin{align*}
q_0(0) = g_0(0), \qquad
q_0(L) = h_0(0).
\end{align*}

\section{A non-effective characterization} \label{noneffectiveapp}
\renewcommand{\theequation}{C.\arabic{equation}}\nequation
The representations of $g_1$ and $h_1$ derived in theorem \ref{th1} yield an effective characterization of the spectral functions $\{A(k), B(k), \mathcal{A}(k), \mathcal{B}(k)\}$ in the sense that the resulting nonlinear system can be solved uniquely at each step in a well-defined perturbative scheme. 

In what follows we present a set of {\it non-effective} formulas (cf. Eqs. (14)-(16)  in  \cite{DMS2001}). Let 
\begin{align}\label{Rcd}
R(0,t,k) = \begin{pmatrix} \overline{d(t, \bar{k})} & c(t, k) \\
\lambda \overline{c(t, \bar{k})} & d(t, k) \end{pmatrix},
\end{align}
where $R$ is defined in (\ref{Rdef}). 
It follows from (\ref{Rmu3relation}) that $R$ satisfies
$$R_t = -2ik^2[\sigma_3, R] + \tilde{Q}(0,t,k)R - Re^{ikL\hat{\sigma}_3} \tilde{Q}(L,t,k).$$
This implies that $\{c, d\}$ satisfy the system
%The function S:= R(0,t,k) = \mu_3(0,t,k)e^{ikL\hat{\sigma}_3}\mu_3(L,t,k)^{-1} satisfies S_t = -2ik^2[\sigma_3, S] + \tilde{Q}(0,t)S - Se^{ikL\hat{\sigma}_3}\tilde{Q}(L,t).
\begin{subequations}\label{alphabetaeqs}
\begin{align} 
&  c_t = -4 ik^2 c -\lambda i(|g_0|^2 + |h_0|^2)c + (2kg_0 + ig_1)d - (2kh_0 + ih_1)e^{2ikL}\bar{d}, 
	\\
&   d_t = \lambda i (|g_0|^2 - |h_0|^2)d + \lambda(2k\bar{g}_0 - i \bar{g}_1)c - \lambda (2kh_0 + ih_1)e^{2ikL}\bar{c}, 
	\\
& c(0,k) = b(k), \qquad d(0,k) = a(k),
\end{align}
\end{subequations}
where $\bar{c}$ and $\bar{d}$ are short-hand notations for $\overline{c(t,\bar{k})}$ and $\overline{d(t,\bar{k})}$, respectively.
Let us consider the Neumann problem. Integration by parts in (\ref{Rdef}) shows that\footnote{In appendix \ref{clemmaapp} we made no assumption on the form of $c(t,k)$ except that it be an entire function satisfying (\ref{cO}). Here, since we are using the form of the function $c(t,k)$ as defined in (\ref{Rcd}), it is much easier to determine its asymptotics.}
$$ c(t, k) = \begin{cases}
 \frac{g_0}{2ik} + O(k^{-2}) + O(e^{2ikL}k^{-1}), & k \in \C^+,  \\
e^{2ikL}\bigl(-\frac{h_0}{2ik} + O (k^{-2})\bigr) + O(k^{-1}), & k \in \C^-,
\end{cases} \qquad k \to \infty.$$
This gives the representations
\begin{subequations}\label{appg0h0rep}
\begin{align}
& g_0(t) = -\frac{1}{\pi} \int_{-\infty}^\infty c_+(t,k) dk,
	\\
& h_0(t) 
%= -\frac{1}{\pi} \int_{-\infty}^\infty [e^{-2ikL}c(t,k) + e^{2ikL}c(t,-k)] dk,$$
= -\frac{1}{\pi} \int_{-\infty}^\infty [\cos(2kL)c_+(t,k) - i\sin(2kL)c_-(t,k)] dk,
\end{align}
\end{subequations}
which can be used to eliminate $\{g_0, h_0\}$ from (\ref{alphabetaeqs}). The resulting system for $\{c, d\}$ is formulated only in terms of $\{g_1, h_1\}$. However, this system is {\it not} effective. Indeed, substituting into (\ref{alphabetaeqs})-(\ref{appg0h0rep}) the expansions
$$d = 1 + \epsilon^2 d_2 + \cdots, \qquad
c = \epsilon c_1 + \epsilon^2Êc_2 + \cdots,$$
the terms of $O(\epsilon^n)$ yield
\begin{align}\label{e4ik2betat}
(e^{4ik^2} c_n)_t = e^{4ik^2}\bigl[2kg_{0n} + ig_{1n} - (2kh_{0n} + ih_{1n})e^{2ikL}\bigr] + \lot
\end{align}
and
$$g_{0n} = -\frac{1}{\pi} \int_{-\infty}^\infty c_{n+} dk, \qquad
h_{0n} = -\frac{1}{\pi} \int_{-\infty}^\infty [\cos(2kL)c_{n+} - i\sin(2kL)c_{n-}] dk.$$
Thus, in order to determine $\{g_{0n}, h_{0n}\}$ at each step of the perturbative scheme, we need to know the combinations $c_{n+}$ and $\cos(2kL)c_{n+} - i\sin(2kL)c_{n-}$. Equation (\ref{e4ik2betat}) shows that these combinations satisfy
\begin{subequations}
\begin{align} \label{combinationsa}
& \bigl\{e^{4ik^2} c_{n+}\bigr\}_t = e^{4ik^2}\bigl[2ig_{1n} - 2kh_{0n}\Delta(k) - ih_{1n}\Sigma(k)\bigr] + \lot,
	\\ \nonumber
& \Bigl\{e^{4ik^2}\bigl[\cos(2kL)c_{n+} - i\sin(2kL)c_{n-}\bigr]\Bigr\}_t 
= e^{4ik^2}\bigl[-2kg_{0n}\Delta(k) + ig_{1n}\Sigma(k) - 2ih_{1n}\bigr] 
	\\ \label{combinationsb}
& \hspace{7.5cm}+ \lot.
\end{align}
\end{subequations}
Although the function $g_{0n}$ has been eliminated from the rhs of (\ref{combinationsa}), the function $h_{0n}$ remains unknown. Similarly, although the function $h_{0n}$ has been eliminated from the rhs of (\ref{combinationsb}), the function $g_{0n}$ remains unknown. 
This shows that the solution is not effective.

\bigskip
\noindent
{\bf Acknowledgement} {\it The authors acknowledge support from the EPSRC, UK. ASF acknowledges support from the Guggenheim foundation, USA.}

\bibliographystyle{plain}
\bibliography{is}

\end{document}